\documentclass{sl}     
\usepackage{slsec}     
\usepackage{stud_log} 
\usepackage{slfoot_arx}
\usepackage{slthm}     
                       
\usepackage[latin9]{inputenc}
\usepackage{longtable}
\usepackage{float}
\usepackage{pmboxdraw}
\usepackage{mathrsfs}
\usepackage{amsmath}
\usepackage{amssymb}
\usepackage{graphicx}

\usepackage{enumerate}

\makeatletter

\providecommand{\tabularnewline}{\\}

\usepackage{enumitem}		
\newlength{\lyxlabelwidth}      

\newenvironment{elabeling}[2][]%
{\settowidth{\lyxlabelwidth}{#2}
\begin{description}[font=\normalfont,style=sameline,
leftmargin=\lyxlabelwidth,#1]}
{\end{description}}

\newcommand{\cml}{[\hspace{-1.5pt}[}
\newcommand{\cmr}{]\hspace{-1.5pt}]}
\usepackage{enumitem}

\AtBeginDocument{
  
}

\theoremstyle{definition}

\HeadingsInfo{Adrian Soncodi}{Automorphisms of the Lattice of Classical Modal Logics}

\begin{document}

\setcounter{page}{1}     

 

\AuthorTitle{Adrian Soncodi}{Automorphisms of the Lattice 
\newline of Classical Modal Logics}

\begin{abstract}
In this paper we analyze the propositional extensions of the minimal classical modal logic system $\mathbf{E}$, which form a lattice denoted as $\mathrm{CExt}\mathbf{E}$.
Our method of analysis uses algebraic calculations with canonical forms, which are a generalization of the normal forms applicable to normal modal logics.
As an application, we identify a group of automorphisms of $\mathrm{CExt}\mathbf{E}$ that is isomorphic to the symmetric group $\mathrm{S}_{4}$.
\end{abstract}

\Keywords{Modal logic, Classical systems, Lattice automorphisms, Uniform replacements.}

\section{Introduction}

The minimal classical modal logic system $\mathbf{E}$ is defined,
for example, in \cite{Chellas}. It extends the classical \emph{propositional
calculus} (PC), from which it inherits the rules of \emph{modus ponens}
(MP) and \emph{uniform substitution} (US). It also adds the \emph{congruence rule} (RE), stated in one of the following equivalent ways:
\begin{elabeling}{00.0000.0000}
\item [{\textsc{(RE-$\lozenge$)}}] \noindent From $\varphi\leftrightarrow\psi$
infer $\lozenge\varphi\leftrightarrow\lozenge\psi$.
\item [{\textsc{(RE-$\square$)}}] From $\varphi\leftrightarrow\psi$ infer
$\square\varphi\leftrightarrow\square\psi$.
\end{elabeling}
This minimal system can in turn be extended by adding to it any set
of well-formed modal formulas as axioms. The resulting classical modal
logic system is the set of formulas that includes all the PC tautologies
as well as the set of axioms and that is closed under MP, US and RE.

In this paper we shall consider only unimodal logics. By analogy
with $\mathrm{NExt}\mathbf{K}$, we denote by $\mathrm{CExt}\mathbf{E}$
the lattice of classical extensions of $\mathbf{E}$. Our method of
analysis is based on the algebra of \textit{canonical forms} for $\mathbf{E}$.
These forms generalize the normal forms described in \cite{Fine}
for $\mathbf{K}$.

In Section 2 we derive the basic properties of canonical forms
and we introduce the concepts of modal context and characteristic minmatrix.
In Section 3 we analyze the effect of a certain group of uniform substitutions
on $\mathbf{E}$-formulas.
We define the concept of prime orbits of canonical forms
and present some related properties.

\pagebreak{}

The main results are in Sections 4 and 5, where we introduce a set of formula
transformations that we call \textit{uniform replacements} and we show that
they determine a group of automorphisms of $\mathrm{CExt}\mathbf{E}$.
Since these automorphisms typically do not preserve normality, the
corresponding lattice symmetries are in fact obscured in $\mathrm{NExt}\mathbf{K}$. 

The method of analysis using normal forms, while not new, is used
quite rarely. In \cite{Moss} the author says that
``\emph{[Kit] Fine's claim that `Normal forms have been comparatively neglected
in the study of modal sentential logic' seems even more cogent thirty
years after its publication}''.
\cite{Moss} actually uses computer calculations with normal forms
to build models for specific systems
like $\mathbf{S4}$, $\mathbf{S4.1}$, $\mathbf{KL}$ etc.~along
the lines presented in Fine's classical paper. Then a number
of properties are derived from these models.

Thus, the normal forms method turns out to be quite fruitful,
and not only for concrete examples but also symbolically,
to prove generic properties.
In our paper we show how symbolic calculations with canonical forms in $\mathbf{E}$
can be used to derive strong results.
Additionally, some of the techniques that we present
can potentially be applied to other logics
(or fragments thereof),
provided that they are algebraizable and have the rule of US.

Interest in the structure of the lattice $\mathrm{NExt}\mathbf{K}$
is shown in \cite{Chagrov1} as well as in its follow-up \cite{Chagrov2}
and in other papers such as \cite{Kracht}.
The the author of the latter uses automorphisms
to investigate $\mathrm{NExt}\mathbf{K4.2}$,
$\mathrm{NExt}\mathbf{K4.3}$ and other sublattices.
At the end he makes several conjectures,
among which one that is considered important and difficult,
namely that the lattice $\mathrm{NExt}\mathbf{K}$
is rigid, i.e. $|Aut(\mathbf{K})|=1$. In this paper we examine $\mathrm{CExt}\mathbf{E}$
and we prove that it is \emph{not} rigid, as it has a non-trivial
group of automorphisms.

We also hope to show here that $\mathrm{CExt}\mathbf{E}$ is
an interesting structure by itself,
and we present some new techniques for studying it.

\section{Normal and Canonical Forms}

\paragraph*{Notations.}

For formulas we use the Greek letters $\varphi,\psi\ldots$ as well as
specific notations detailed further on. Propositional variables are
$p,q,r\ldots$ and the logical constants are \emph{0} and \emph{1}
(note the italics).
The operators are listed in Table 1, in descending order of precedence.

Whenever necessary, parentheses are used for disambiguation. Examples
of \emph{well-formed formulas} (wff) are $p!q+qr\rightarrow\mbox{!}qr$
and $(\square p\leftrightarrow\square\lozenge p)\rightarrow\lozenge\square p$.

The \emph{modal degree} of a formula is the largest number of nested
modal operators found in any sub-formula of the given formula. We
also refer to it as the \emph{level} of the formula. In the above
examples, the modal degrees (levels) are 0 and 2 respectively.

\pagebreak{}

\begin{center}
\begin{tabular}{|c|l|}
\noalign{\vskip1pt}
\multicolumn{2}{c}{{\footnotesize{Table 1: Operators precedence}}}\tabularnewline[1pt]
\noalign{\vskip1pt}
\hline 
\noalign{\vskip1pt}
{\scriptsize{$\square$, $\lozenge$ and $!$}} & \emph{\scriptsize{necessity}}{\scriptsize{, }}\emph{\scriptsize{possibility}}{\scriptsize{
and }}\emph{\scriptsize{negation}}\tabularnewline[1pt]
\noalign{\vskip1pt}
\hline 
\noalign{\vskip1pt}
{\scriptsize{$\circ$ and $*$}} & \emph{\scriptsize{uniform substitution (US)}}{\scriptsize{ and }}\emph{\scriptsize{uniform
replacement (UR) }}{\scriptsize{application}}\tabularnewline[1pt]
\noalign{\vskip1pt}
\hline 
\noalign{\vskip1pt}
{\scriptsize{$\cdot$}} & \emph{\scriptsize{conjunction}}{\scriptsize{, as well as}}\emph{\scriptsize{
US}}{\scriptsize{ and }}\emph{\scriptsize{UR composition}}{\scriptsize{
(typically omitted though)}}\tabularnewline[1pt]
\noalign{\vskip1pt}
\hline 
\noalign{\vskip1pt}
{\scriptsize{$+$}} & \emph{\scriptsize{disjunction}}\tabularnewline[1pt]
\noalign{\vskip1pt}
\hline 
\noalign{\vskip1pt}
{\scriptsize{$\rightarrow$ and $\leftrightarrow$}} & \emph{\scriptsize{logical implication}}{\scriptsize{ and }}\emph{\scriptsize{logical
equivalence}}{\scriptsize{ within formulas}}\tabularnewline[1pt]
\noalign{\vskip1pt}
\hline 
\noalign{\vskip1pt}
{\scriptsize{$=$ and $\approx_{_\mathbf{B}}$}} & \emph{\scriptsize{identity}}{\scriptsize{ and}}\emph{\scriptsize{
equivalence (equiprovability in $\mathbf{B}$)}}{\scriptsize{ of formulas;
$\mathbf{B}$ omitted when $\mathbf{E}$}}\tabularnewline[1pt]
\hline 
\noalign{\vskip1pt}
\end{tabular}\vspace{-0.1mm}

\par\end{center}
 
Let $\mathcal{F}(v,d)$ be set of all unimodal wffs in a number of
variables not exceeding $v$ and of modal degree not exceeding $d$.

Regarding the above notations as well as those introduced further on,
we wish to mention here that the notations for the
present paper have been a challenge and are thus the result of an
extensive consideration. The problem is that in our case the traditional notations
often result in very long formulas. Since our paper consists mainly of algebraic
calculations with normal or canonical forms, which are essentially
sums of products, we opted for a more compact representation. Therefore
we kindly ask the reader to bear with us in terms of these notations,
because we believe that in the end it is worth it,
for the sake of compactness and readability of the proofs.

\paragraph*{Disjunctive canonical form.}

The modal \emph{disjunctive normal form} (DNF) for normal modal logic formulas
is described in \cite{Fine}.
This DNF representation of a formula in $\mathcal{F}(v,d)$ is unique up to the reordering
allowed by the commutative logical connectors, and it is also equiprovable
to the initial formula. An example $\mathcal{F}(1,1)$ formula in modal DNF
is as follows:\vspace{-3mm}

\[
\varphi(p)=p\lozenge p\lozenge!p+p\lozenge p!\lozenge!p+\mbox{!}p\lozenge p\lozenge!p+\mbox{!}p!\lozenge p\lozenge!p
\]

The modal DNF is a sum (disjunction) of \emph{minterms}.
(Note that in some papers minterms are called ``full normal forms'';
but we reserve the phrase ``normal forms'' to
denote the generic DNF representation of formulas.)
Each minterm is a product (conjunction) of  modal and non-modal \emph{factors}.
The above formula can also be represented
in a matrix format as follows:\vspace{-2mm}

\begin{equation}
\varphi(p)=\begin{array}{c|cccc|}
p & 1 & 1 & 0 & 0\\
\hline \lozenge p & 1 & 1 & 1 & 0\\
\lozenge!p & 1 & 0 & 1 & 1
\end{array}\label{eq:CMM-T}
\end{equation}

Because of this representation, we also call the normalized formula
a \emph{minmatrix}. The binary columns correspond to the formula's
minterms. Left of the binary entries, the labels denote the modal
and non-modal factors, below and above the line respectively.
A minmatrix entry is the factor's \emph{state} in the corresponding minterm: 0
if the row's factor occurs complemented, otherwise 1. The product
of all the non-modal factors from a minterm, each in their respective
state, constitutes its (level 0) \emph{prefix},
which is in fact a Boolean DNF minterm. Barring the reordering of rows and columns,
in every context the minmatrix representation of a formula is unique.

However, we need to point out that the above matrix is only a pictorial
representation of a normalized formula. It is useful to quickly show
what is relevant, namely the states of the minterm factors, but it
can be equivalently replaced by the algebraic formula at any time.
Other than that, a minmatrix \emph{is} in fact a normalized formula. 

We now introduce the \emph{disjunctive canonical form} (DCF) for  $\mathbf{E}$.
Fix the number of variables $v$.
For $\mathcal{F}(v,0)$, the DCF formulas are the Boolean DNF formulas. 
When $d>0$, a DCF formula from $\mathcal{F}(v,d)$ is a similar sum of minterms
with various state combinations for the factors, which now include:\vspace{-1mm}
\begin{itemize}
\item All \emph{non-modal factors} consisting of all the $v$ propositional
variables.\vspace{-0.5mm}
\item All \emph{modal factors} $\lozenge\phi$, where $\phi$ is
every DCF \emph{formula} from $\mathcal{F}(v,d-1)$.
(Note the difference from the modal DNF, where the modal factors are $\lozenge\mu$
for every \emph{minterm} $\mu$ from $\mathcal{F}(v,d-1)$.)
\vspace{-1mm}
\end{itemize}
To convert an $\mathbf{E}$-formula to an equiprovable DCF,
one uses RE and the rule of \emph{substitution of equivalents}
(EQ, also called REP in of \cite{Chellas}, where it is proven for $\mathbf{E}$ in Theorem 8.3).
Given a level $d$ formula $\varphi$ with $d>0$,
the recursive DCF conversion procedure consists of the following steps:\vspace{-1mm}
\begin{enumerate}
\item For all outermost modalities $\lozenge\phi$,
convert $\phi$ to the level $d-1$ DCF.\vspace{-0.5mm}
\item Apply the Boolean DNF conversion procedure to the resulting formula,
for the purpose of which propositional variables and distinct outermost modalities
are considered atomic formulas and any missing level $d$ modal factors
$\lozenge\theta$ are re-introduced using EQ with $\varphi\approx\varphi(\lozenge\theta+\mbox{!}\lozenge\theta)$.
\vspace{-1mm}
\end{enumerate}
Unlike the DNF conversion from $\mathbf{K}$, since $\mathbf{E}$ lacks axioms N and K,
this procedure does \emph{not} (and cannot) use
 $\lozenge(\phi+\theta)\approx\lozenge\phi+\lozenge\theta$ to split modalities,
and it must also preserve $\lozenge\mathit{0}$,
which is not equiprovable to $\mathit{0}$ in $\mathbf{E}$.

Below we show an $\mathcal{F}(1,1)$ formula in DCF and its minmatrix form:\vspace{-1mm}
\[
\varphi(p)=p\lozenge\mathit{1}\lozenge p\lozenge!p\lozenge\mathit{0}+p!\lozenge\mathit{1}\lozenge p!\lozenge!p!\lozenge\mathit{0}+\mbox{!}p!\lozenge\mathit{1}!\lozenge p\lozenge!p\lozenge\mathit{0}+\mbox{!}p!\lozenge\mathit{1}\lozenge p\lozenge!p!\lozenge\mathit{0}
\]
\vspace{-3mm}
\[
\varphi(p)=\begin{array}{c|cccc|}
p & 1 & 1 & 0 & 0\\
\hline \lozenge\mathit{1} & 1 & 0 & 0 & 0\\
\lozenge p & 1 & 1 & 0 & 1\\
\lozenge!p & 1 & 0 & 1 & 1\\
\lozenge\mathit{0} & 1 & 0 & 1 & 0
\end{array}
\]

Note how the DCF conversion procedure \emph{promotes} DCF formulas from levels 0, 1, ... $d-1$ to $d$ as needed. Hence the set of level $d+1$ modal factors
includes equiprovables of all the modal factors from levels 0 to $d$.
Then by reordering the factors
we define a \emph{level} $i$ \emph{minterm prefix}, $0\le i\le d$,
as the product of all the factors that are equiprovable to formulas of levels $\le i$.

By PC, every level $d$ minterm $\mu$ is equiprovable to a sum of level $d+1$ minterms,
which have the same level $d$ prefix
and all the state combinations for the remaining level $d+1$ modal factors.
We call them the \emph{immediate descendants} of $\mu$,
and $\mu$ their \emph{immediate ancestor}.
Since distinct level $d$ minterms disagree on at least one of their factors' states,
they have disjoint sets of immediate descendants.
Then we can define descendence-ancestry between any levels,
and the minterms from all $\mathcal{F}(v,i)$,
$0\le i \le d$, can be represented as $2^v$ ancestry trees rooted in the level 0 minterms.
\vspace{1mm}

\begin{elabeling}{00.0000.0000.00}
\item [{\textbf{\textsc{Theorem~1}}}] \noindent
\emph{A DCF formula from $\mathcal{F}(v,d)$ is a theorem of $\mathbf{E}$
iff it is the sum of all the DCF minterms from $\mathcal{F}(v,d)$.}
\end{elabeling}\vspace{-1mm}
\textbf{\textsc{Proof.}}\qquad{}
Sufficiency is a consequence of the fact that the sum of all the minterms from $\mathcal{F}(v,d)$ is a Boolean tautology.
For the necessity we use the fact that by Theorem 9.8 of \cite{Chellas},
$\mathbf{E}$ is complete with respect to the class of all neighborhood frames.
If a theorem did not include all the minterms,
then the remaining minterms would be unsatisfiable formulas.
So it suffices to show that every DCF minterm is valid at some world in some model.

For this, we construct models similar to the graded models in \cite{Fine},
but we adapt them to neighborhood frames.
Fix $v$ and $d$. We take the set of worlds $W=\bigcup_{0\leq i\leq d}W^{i}$,
where $W^{i}=\{w_{\mu}:\mu\text{ is a DCF minterm from }\mathcal{F}(v,i)$\}.
Worlds can then be mapped to the minterms in the ancestry trees defined
above, and we adopt a similar terminology for them. Let $F=(W,N)$
be a neighborhood frame, with $N:W\rightarrow\wp(\wp(W))$
its neighborhood function, i.e.~$N(w_{\mu})$ contains the neighborhoods
of $w_{\mu}$. In a model $\mathscr{M}=(F,V)$ based on this frame,
the valuation function $V:W\times\mathcal{F}(v,d)\rightarrow\{\mathit{0},\mathit{1}\}$
is defined in the usual way and
$V(w_{\mu},\lozenge\phi)=\mathit{1}$ iff
$W\setminus\{w\in W:V(w,\phi)=\mathit{1}\}\notin N(w_{\mu})$.

At each $w_{\mu}$ we take $V(w_{\mu},p_{i})$ to match the state
of the variable $p_{i}$ in the level 0 prefix of $\mu$, i.e.~precisely
the valuations that make the level 0 prefix of $\mu$ valid at $w_{\mu}$.
The goal is to define $\mathscr{M}$ such that $V(w_{\mu},\mu)=\mathit{1}$
for all $w_{\mu}$. But we note that in this case $\mathscr{M}$ must meet more stringent
requirements:\vspace{-0.5mm}
\begin{itemize}
\item If $V(w_{\mu},\mu)=\mathit{1}$ for some $w_{\mu}$,
then $V(w_{\nu},\mu)=\mathit{1}$ for all descendants $\nu$ of $\mu$.
This is because $\mu$ is (equiprovable to) a prefix of $\nu$.\vspace{-0.5mm}
\item If $\mu$ and $\mu'$ are distinct level $i$ minterms, $0\leq i\leq d$,
then they disagree on at least one of the states of their factors,
hence $V(w_{\mu},\mu')=\mathit{0}$ and $V(w_{\nu},\mu')=\mathit{0}$
for all descendants $\nu$ of $\mu$.
\item If $\mu$ is a level $j$ minterm and $j<i$,
then only one level $i$ descendant $\nu$ of $\mu$ has $V(w_\mu,\nu)=\mathit{1}$.
Here we have a choice, but it is convenient to require $\mathscr{M}$ to be such that
this is always the special minterm $\nu$
whose states of all the modal factors occurring \emph{after} its level $j$ prefix are 0.
\end{itemize}\vspace{-1mm}
For $\varphi\in\mathcal{F}(v,d)$,
let $X(\varphi)=\{w\in W:V(w,\varphi)=\mathit{1}\}$.
Then for any minterm $\mu$, the above requirements uniquely determine $X(\mu)$,
hence if $\phi=\Sigma\mu_k$, $X(\phi)=\cup X(\mu_k)$.
And by examining $X(\mu)$ at all the worlds in the ancestry trees of $\mathscr{M}$
we see that, recursively (by level),
the above requirements also have the following consequences:\vspace{-1mm}

\begin{elabeling}{0000.}
\item [{C1:}] \noindent If $\phi$ and $\phi'$ are distinct level
$i$ DCF formulas, then $X(\phi)\not=X(\phi')$, since they include
sets of worlds rooted at different level $i$ worlds.
\item [{C2:}] If $\theta$ is a level $j$ DCF formula, $j<i$, then $X(\phi)=X(\theta)$
iff $\phi$ is the sum of all level $i$ descendants of the minterms
of $\theta$, therefore $\phi\approx\theta$.
\end{elabeling}\vspace{-1.5mm}
We now show that we can choose $N(w)$ at all $w$ such that
$\mathscr{M}$ meets our requirements.
For each $\mu$, to assign $N(w_{\mu})$ we start from $\wp(W)$
and we eliminate neighborhoods as needed to make the modal factors valid at $w_\mu$. 

For a level 0 $\mu$, already $V(w_{\mu},\mu)=\mathit{1}$,
so we take $N(w_{\mu})=\wp(W)$.
Then $V(w_{\mu},\lozenge\phi)=\mathit{0}$ for all $\phi$,
so only suitable special minterms are valid at $w_{\mu}$.

For a level 1 $\mu$, again its level 0 prefix is already valid at $w_{\mu}$.
For each of its modal factors $\lozenge\phi$,
if its state in $\mu$ is 1 then we exclude neighborhood $W\setminus X(\phi)$ from $N(w_\mu)$,
otherwise we don't.
We do the same at all $w_\nu$ for all descendants $\nu$ of $\mu$.
Since C1 holds on level 0, there is no conflict between these inclusions/exclusions.
And for a special level 1 minterm we exclude nothing from $\wp(W)$,
so this is valid at $w_{\mu}$ and at its ancestor world.
Then after processing all $\mu$ on level 1,
their $X(\mu)$ are as per our target $\mathscr{M}$.

For the induction step we assume that for all $\mu$ up to level $i$,
$X(\mu)$ is as per our target $\mathscr{M}$.
Consider a level $i+1$ minterm $\mu$.
Its level $i$ prefix is already valid at $w_\mu$,
so we exclude or include neighborhoods at $w_\mu$ and at all its descendants
as per the states of the \emph{remaining} level $i+1$ modal factors $\lozenge\phi$ of $\mu$.
If $\mu$ is a special minterm then there are no new sets to exclude,
hence $\mu$ is valid at the immediate ancestor world.
Otherwise, since C1 and C2 hold on level $i$,
there is no conflict on level $i+1$ with neighborhoods already included/excluded.
(From C2, if $X(\phi)=X(\theta)$ for a lower level $\theta$,
then $\phi\approx\theta$ and $\lozenge\phi$ is already part of the level $i$ prefix of $\mu$.)
Then after processing all $\mu$ on level $i+1$,
their $X(\mu)$ are again as per our target $\mathscr{M}$.

We iterate this processing up to level $d$,
and so by construction $\mathscr{M}$ meets our goal.
Hence every minterm is a satisfiable formula.
\textsc{\scriptsize{\hfill{}$\blacksquare$}}\vspace{2mm}

As an immediate consequence, distinct level $d$ DCF formulas cannot be equiprovable
(otherwise there would be minterms $\mu\approx\mathit{0}$).

\pagebreak{}

\paragraph*{Modal contexts.}

We denote systems by boldface and axioms by Roman letters. $\mathbf{F}$
is reserved for the inconsistent system (consisting of all formulas).

Consider a modal logic system $\mathbf{B}$ as a basis for the discussion.
In our case, $\mathbf{B}$ is either $\mathbf{K}$ or $\mathbf{E}$.
We define a \emph{modal context} $\mathbf{B}[v,d]$ as the quotient
$\mathcal{F}(v,d)/\approx_{\hspace{-1pt}_{\mathbf{B}}}$, i.e.~the Lindenbaum-Tarski
algebra of classes of $\mathbf{B}$-equiprovable formulas. But we
tacitly equate a class with a representative from it, so that we can
still refer to these classes as formulas (or minmatrices).
This is similar to writing
$\mbox{1+1=0}$ instead of $\hat{1}+\hat{1}=\hat{2}=\hat{0}$ in $\mathbb{Z}_{2}$,
which is often done and has the advantage of avoiding hats over large formulas.

Since there is only one DCF representative per class,
when we want to emphasize that we refer to it
we write $[\varphi]$ (or $_{v}^{d}[\varphi]$ to specify the context).
Yet we always have $\varphi\approx[\varphi]$.
The formula $\mathit{0}$ corresponds to the empty minmatrix $[\mathit{0}]$ and the formula containing all the minterms
from the context by the minmatrix $[\mathit{1}]$
(notation not to be confused with the reference \cite{Blackburn}).

With this convention,
$\mathbf{E}[v,d]$ is the (finite) set of all minmatrices from the context.
Then the Boolean operations on formulas can be performed as set operations
on the corresponding minmatrix minterms;
namely, union for disjunction, intersection for conjunction, complementation
(with respect to $[\mathit{1}]$) for negation etc.

For this reason we may also interpret $[\varphi]$ as a set (rather than a sum) of minterms.
This notation overloading allows us to avoid the constant use of conversion operators
between DCF formulas and their sets of minterms,
while it can still be disambiguated from the surrounding text
(for example, in $f:[\mathit{1}]\rightarrow[\mathit{1}])$.
And it allows us to be brief by writing ``minmatrix intersection'' instead of
``the minmatrix that is the sum of the minterms
from the intersection of the sets of minterms from ...''

With this notation we can also write, for example, $[\varphi]\subset[\psi]$.
This partial order relationship determines a lattice structure on $\mathbf{E}[v,d]$.
\vspace{-2mm}

\paragraph*{Characteristic minmatrix and $\mathbf{E}[v,d]$ systems.}

From Theorem 1:\vspace{-3mm}

\begin{equation}
\vdash\mathbf{_{\hspace{-3pt}_{E}}\,}\varphi\rightarrow\psi\quad\text{iff}\quad[\varphi\rightarrow\psi]=[\mathit{1}]\quad\text{iff}\quad[\varphi]\subseteq[\psi]\label{eq:basic-1}
\end{equation}\vspace{-6mm}
\begin{equation}
\vdash\mathbf{_{\hspace{-3pt}_{E}}}\,\varphi\leftrightarrow\psi\quad\text{iff}\quad[\varphi\leftrightarrow\psi]=[\mathit{1}]\quad\text{iff}\quad[\varphi]=[\psi]\label{eq:basic-2}
\end{equation}
\vspace{-6mm}

Let $\mathbf{S\in\mathrm{CExt}\mathbf{E}}$. The\emph{ characteristic
minmatrix} (CMM) of $\mathbf{S}$ for a context $\mathbf{E}[v,d]$,
denoted as $\cml\mathbf{S}\cmr$ or $_{v}^{d}\cml\mathbf{S}\cmr$,
is defined as the minmatrix intersection (conjunction)
of all the $\mathbf{S}$-theorems from the context.
The definition is sound,
since there is only a finite number of equiprovable $\mathbf{S}$-theorems per context.

Note that the minmatrices used in the intersection must belong to the context,
but the formal proof of the corresponding $\mathbf{S}$-theorems may
involve formulas from other contexts (and in fact this is often necessary).

This definition and (\ref{eq:basic-1}) imply that
$_{v}^{d}\varphi$ is a theorem of $\mathbf{S}$ iff $\,{}_{v}^{d}\cml\mathbf{S}\cmr\subseteq{}_{v}^{d}[\varphi]$.
When this is the case, one can construct a formal proof for $\varphi$
from $\cml\mathbf{S}\cmr$ using mainly EQ
and the \emph{PC monotony} rule (from $\vdash p$ infer $\vdash p+q$).

A system is determined by the set of its CMMs from all the contexts.
But within a given context, distinct systems may share the same CMM
(when they prove the same theorems within that context). Yet we can
always associate a unique system with a minmatrix $[\mathrm{S}]$
that is a CMM. This is the system $\mathbf{S}$ that extends $\mathbf{E}$
precisely with $[\mathrm{S}]$ (i.e.~axiom $\mathrm{S}$), and we
call it \emph{the} $\mathbf{E}[v,d]$ \emph{system} corresponding to $[\mathrm{S}]$.
Conversely, for an $\mathbf{E}[v,d]$ system $\mathbf{S}$,
we denote by $[\mathrm{S}]$ or $\cml\mathrm{S}\cmr$ its determining
axiom (or CMM) from that context.

On the other hand, not every minmatrix $[\mathrm{S}]\in\mathbf{E}[v,d]$
can be a CMM. If we can derive from $\mathrm{S}$
(e.g.~by US) another formula $\varphi$
such that $[\mathrm{S\,}\varphi]\subsetneq[\mathrm{S}]$, then $[\mathrm{S}]$
cannot be $\cml\mathrm{\mathbf{S}}\cmr$ for any system $\mathbf{S}$, since by definition
$\cml\mathrm{\mathbf{S}}\cmr\subseteq[\mathrm{S}\varphi]$.

Let $\mathbf{E}\cml v,d\cmr$ be the set of CMMs from $\mathbf{E}[v,d]$,
where $\mathbf{E}\cml v,d\cmr\subset\mathbf{E}[v,d]$.
The partial order induced by the set inclusion relationship between CMMs
determines a lattice structure on $\mathbf{E}\cml v,d\cmr$,
where in every context $\cml\mathbf{F}\cmr=[\mathit{0}]$
and $\cml\mathbf{E}\cmr=[\mathit{1}]$, namely~$\mathbf{F}$ and $\mathbf{E}$
are $\mathbf{E}[v,d]$ systems, as well as 
the $\bot$ and $\top$ elements of the $\mathbf{E}\cml v,d\cmr$ lattice respectively.
But we need to point out the difference between
the set-based operations in $\mathbf{E}[v,d]$
and the lattice operations in $\mathbf{E}\cml v,d\cmr$.
Denote the latter by $\vee$ and $\wedge$.

\bgroup \renewcommand*\theenumi{\alph{enumi}}\renewcommand*\labelenumi{\theenumi)}\vspace{1mm}

\begin{elabeling}{00.0000.0000.00}
\item [{\textbf{\textsc{Theorem~2}}}] \noindent
\emph{Let }$\cml\mathbf{S'}\cmr$\emph{and }$\cml\mathbf{S''}\cmr$
\emph{ be CMMs from }$\mathbf{E}\cml v,d\cmr$\emph{. Then:}
\end{elabeling}\vspace{-3mm}

\begin{enumerate}
\item \emph{$\cml\mathbf{S'}\cmr\vee\cml\mathbf{S''}\cmr=\cml\mathbf{S'}\cmr\cup\cml\mathbf{S''}\cmr$}

\item \emph{$\cml\mathbf{S'}\cmr\wedge\cml\mathbf{S''}\cmr\subseteq\cml\mathbf{S'}\cmr\cap\cml\mathbf{S''}\cmr$}
\end{enumerate}
\textbf{\textsc{Proof.}}\qquad{}
Denote by $[\mathrm{S'}]\triangleq\cml\mathbf{S'}\cmr$
and $[\mathrm{S}'']\triangleq\cml\mathbf{S''}\cmr$ the determining
axioms of $\mathbf{S'}$ and $\mathbf{S''}$.
For a), the inclusion $\supseteq$ is obvious,
so we must show that $[\mathrm{S}]\triangleq\cml\mathbf{S'}\cmr\cup\cml\mathbf{S''}\cmr$
is indeed a CMM.
Let $\mathbf{S}$ be the $\mathbf{E}[v,d]$ system
corresponding to $[\mathrm{S}]$. If $[\mathrm{S}]$ is not a CMM,
then there is a theorem $\varphi$ of $\mathbf{S}$ such that $[\mathrm{S\,\varphi}]\subsetneq[\mathrm{S}]$.
Then at least one of $[\mathrm{S'\,\varphi}]\subsetneq[\mathrm{S'}]$
or $[\mathrm{S''\,\varphi}]\subsetneq[\mathrm{S''}]$ holds.
But by (\ref{eq:basic-1}) both $\mathbf{S'}$ and $\mathbf{S''}$
already prove $\mathrm{S}$, hence any theorem of $\mathbf{S}$, including
$\varphi$, so at
least one of $\cml\mathbf{S'}\cmr$ and\emph{ }$\cml\mathbf{S''}\cmr$
is not a CMM, contradicting our assumption.
Next, b) holds because any system that proves $[\mathrm{S'}]$ and $[\mathrm{S}'']$ proves at least $\cml\mathbf{S'}\cmr\cap\cml\mathbf{S''}\cmr$. But here the inclusion may be strict, since the combination of $[\mathrm{S'}]$ and $[\mathrm{S}'']$ may prove a CMM that is stronger than this intersection.
\textsc{\scriptsize{\hfill{}$\blacksquare$}}\vspace{3mm}
\egroup

As defined, an $\mathbf{E}[v,d]$ system is finitely-axiomatizable,
and it is also the weakest extension of $\mathbf{E}$ that has that
CMM in the given context. Obviously, every finitely-axiomatizable
system is an $\mathbf{E}[v,d]$ system in some context(s).
A system that is not finitely-axiomatizable is not an $\mathbf{E}[v,d]$ system
in any context, but it still has a CMM in every context.

\pagebreak

\section{The Minterm Structure of CMMs}

In this section we take a closer look at CMMs to determine what specific sets of minterms
they may consist of.
We derive a necessary condition for a minmatrix to
be the CMM of some system.

\paragraph*{Uniform substitutions.}

In the following we shall assume that the working context $\mathbf{E}[v,d]$
can accommodate all the formulas involved. We write $\langle\alpha_{i}\rangle$
as a shorthand for $(\alpha_{1},\ldots,\alpha_{v})$, for
example $\varphi(p_{1},\ldots,p_{v})=\varphi\langle p_{i}\rangle$
and $(\sigma_{1}\langle p_{j}\rangle,\ldots,\sigma_{v}\langle p_{j}\rangle)=\langle\sigma_{i}\langle p_{j}\rangle\rangle$,
even though not all these formulas necessarily depend on all the propositional
variables $p_{i}$, $1\leq i\leq v$.

Let $\sigma$ be the uniform substitution $\langle p_{i}\rangle\mapsto\langle\sigma_{i}\langle p_{j}\rangle\rangle$,
where $\sigma_{i}$ are formulas, $1\leq i\leq v$. The result, denoted
as $\varphi\circ\sigma$, of applying the substitution $\sigma$ to
a formula $\varphi=\varphi\langle p_{i}\rangle$ is the formula $(\varphi\circ\sigma)\langle p_{i}\rangle$
obtained by consistently replacing every occurrence of every propositional
variable $p_{i}$ in $\varphi$ by the corresponding $\sigma_{i}\langle p_{j}\rangle$.
We write this as $\varphi\circ\sigma=(\varphi\circ\sigma)\langle p_{i}\rangle=\varphi\langle\sigma_{i}\langle p_{j}\rangle\rangle$.
Formally, this operation is defined by the following rules, applied
recursively to the sub-formulas $\psi$, $\theta$, $\ldots$ that
occur in $\varphi$:
\begin{elabeling}{00.00.0000}
\item [{(US-1)}] $\mathit{0}\circ\sigma\triangleq\mathit{0}$ and $\mathit{1}\circ\sigma\triangleq\mathit{1}$
\item [{(US-2)}] $p_{i}\circ\sigma\triangleq\sigma_{i}\langle p_{j}\rangle$
\item [{(US-3)}] $(\psi+\theta)\circ\sigma\triangleq\psi\circ\sigma+\theta\circ\sigma$
\item [{(US-4)}] $(!\psi)\circ\sigma\triangleq\mbox{!}(\psi\circ\sigma)$
\item [{(US-5)}] $(\lozenge\psi)\circ\sigma\triangleq\lozenge(\psi\circ\sigma)$
\end{elabeling}

Then the following are immediate consequences:
\begin{elabeling}{00.00.0000}
\item [{(US-6)}] $(\psi\,\theta)\circ\sigma\approx(\psi\circ\sigma)(\theta\circ\sigma)$
\item [{(US-7)}] $(\psi\rightarrow\theta)\circ\sigma\approx\psi\circ\sigma\rightarrow\theta\circ\sigma$
\item [{(US-8)}] $(\psi\leftrightarrow\theta)\circ\sigma\approx\psi\circ\sigma\leftrightarrow\theta\circ\sigma$
\item [{(US-9)}] $(\psi\not\leftrightarrow\theta)\circ\sigma\approx\psi\circ\sigma\not\leftrightarrow\theta\circ\sigma$
\item [{(US-10)}] $(\square\psi)\circ\sigma\approx\square(\psi\circ\sigma)$
\end{elabeling}

Also, using US-1 to US-10 above and the fact that $\varphi\approx\psi$
is defined as $\vdash\mathbf{_{\hspace{-2pt}_{\mathbf{E}}}\,}\varphi\leftrightarrow\psi$,
one can easily prove the following additional properties:
\begin{elabeling}{00.00.0000}
\item [{(US-11)}] If $\varphi\approx\psi$ then $\varphi\circ\sigma\approx\psi\circ\sigma$
\item [{(US-12)}] If $\varphi\,\psi\approx\mathit{0}$ then $(\varphi\circ\sigma)(\psi\circ\sigma)\approx\mathit{0}$
\end{elabeling}

Our immediate interest is in level 0 substitutions, i.e.~substitutions
where all $\sigma_{i}$ are level 0 formulas. We also call them \emph{context-preserving}
substitutions, because by applying them to any formula $\varphi$,
neither the number of variables nor the modal degree increase.
If $v$ or $d$ actually decrease for a sub-formula of $\varphi$,
it can always be promoted back to an equiprovable formula from $\mathbf{E}[v,d]$
using the tautology $\varphi\approx\varphi(\psi+{!}\psi)$ and EQ.

There are $2^{2^{v}}$ formulas in $\mathbf{E}[v,0]$, therefore
$2^{v\cdot2^{v}}$ context-preserving substitutions $\sigma$ that
can be applied to any $\varphi\in\mathbf{E}[v,d]$. They are defined
independently of the formulas $\varphi$ of various modal degrees.
Let $\mathcal{S}(v,0)$ be the set of all level 0 substitutions
in $v$ variables. The \emph{composition} $\sigma\sigma'$ of substitutions
$\sigma=\langle\sigma_{i}\langle p_{j}\rangle\rangle$
and $\sigma'=\langle\sigma_{i}'\langle p_{j}\rangle\rangle$ from $\mathcal{S}(v,0)$
is defined as follows:\vspace{-3mm}

\[
(\sigma\sigma')\langle p_{i}\rangle\triangleq\langle\sigma_{i}\langle\sigma_{j}'\langle p_{k}\rangle\rangle\rangle
\]

The composition of level 0 substitutions is obviously well-defined
(the result is context-preserving) and its associativity is straightforward
to verify. With this operation $\mathcal{S}(v,0)$ is a monoid, whose
unit is the identical substitution $\varsigma_{0}=\varsigma_{0}\langle p_{i}\rangle\triangleq\langle p_{i}\rangle$.
Then $\varphi\circ\sigma$ actually defines a (right) monoid action
of $\mathcal{S}(v,0)$ on $\mathbf{E}[v,d]$, compatibility being
ensured since for any formula $\varphi\in\mathbf{E}[v,d]$ we have:\vspace{-3mm}

\[
\varphi\circ(\sigma\sigma')=\varphi\langle\sigma_{i}\langle\sigma_{j}'\langle p_{k}\rangle\rangle\rangle=\varphi\langle\sigma_{i}\langle p_{j}\rangle\rangle\circ\sigma'=(\varphi\circ\sigma)\circ\sigma'
\]
\vspace{-8mm}

\paragraph{Prime substitutions and prime orbits.}

Assume that $\mathbf{S}$ extends $\mathbf{E}$ with axiom $\mathrm{S}$.
From $\vdash\mathbf{_{\hspace{-2pt}_{S}}\,}\mathrm{S}$ infer $\vdash\mathbf{_{\hspace{-2pt}_{S}}\,}\varphi$
for some $\varphi$, then $\vdash\mathbf{_{\hspace{-2pt}_{S}}\,}\mathrm{S\,\varphi}$,
with $[\mathrm{S\,\varphi}]=[\mathrm{S}]\cap[\varphi]$. If $[\mathrm{S}\,\varphi]\subsetneq[\mathrm{S}]$,
we say that $[\mathrm{S}]$ \emph{collapses} (by intersection with
the minmatrix of some other theorem). Let $\mathrm{\varphi=S\circ\sigma}$
for a context-preserving $\sigma$. If $[\mathrm{S\,}\varphi]\subsetneq[\mathrm{S}]$
we say that $[\mathrm{S}]$ \emph{collapses under} $\sigma$; otherwise
if $[\mathrm{S}\,\varphi]=[\mathrm{S}]$, i.e.~$[\mathrm{S}]\subseteq[\varphi]$,
we say that $[\mathrm{S}]$ is \emph{immune} to $\sigma$.

As we have seen, if $[\mathrm{S}]$ collapses then it can not be a
CMM, since this requires $\mathrm{[S}]\subseteq[\mathrm{S}\,\varphi]$.
Thus, a candidate CMM must first of all be immune to all the context-preserving
substitutions. This is a necessary, albeit not sufficient condition
for a minmatrix to be a CMM.

We therefore analyze which minmatrices are immune to level 0 substitutions.
We begin by considering a subset of level 0 substitutions that we
call \emph{prime substitutions}. They are defined as the invertible
elements of the monoid $\mathcal{S}(v,0)$, hence they form a group
that we denote by $\mathcal{S}_{p}(v,0)$.

It turns out that prime substitutions are precisely the level 0 substitutions
that always transform a minterm into a single minterm. As such, they
generate automorphisms of the lattice $\mathbf{E}[v,d]$. The
following theorems establish this result.

\pagebreak{}

\begin{elabeling}{00.0000.0000.0000}
\item [{\textbf{\textsc{Theorem~3}}}] \emph{For any given $v$, the group
}$\mathcal{S}_{p}(v,0)$\emph{ is isomorphic to the} \emph{symmetric
group }$\mathrm{S}_{2^{v}}$\emph{.}
\end{elabeling}
\textbf{\textsc{Proof.}}\qquad{}For every $\varsigma\in\mathcal{S}_{p}(v,0)$,
let $f_{\varsigma}:\mathbf{E}[v,0]\rightarrow\mathbf{E}[v,0]$ be
defined as $f_{\varsigma}(\varphi)\triangleq\varphi\circ\varsigma$.
We prove that $f_{\varsigma}$ is an automorphism of the lattice $\mathbf{E}[v,0]$.
First, we show that it is injective. Assume $\varphi,\psi\in\mathbf{E}[v,0]$
and $f_{\varsigma}(\varphi)\approx f_{\varsigma}(\psi)$, i.e.~$\varphi\circ\varsigma\approx\psi\circ\varsigma$.
Since $\varsigma$ is invertible, we have $\varphi\circ\varsigma\circ\varsigma^{-1}\approx\psi\circ\varsigma\circ\varsigma^{-1}$
and by compatibility $\varphi\circ(\varsigma\,\varsigma^{-1})\approx\psi\circ(\varsigma\,\varsigma^{-1})$,
i.e.~$\varphi\approx\psi$. Next, being injective on the finite set
$\mathbf{E}[v,0]$, $f_{\varsigma}$ must be a bijection. Properties
US-1, US-3 and US-6 show that $f_{\varsigma}$ is also compatible
with the lattice operations in $\mathbf{E}[v,0]$, hence it is an
automorphism.

We observe that if the prime substitution $\varsigma$ is $\langle p_{i}\rangle\mapsto\langle\varsigma_{i}\langle p_{j}\rangle\rangle$,
then we have $f_{\varsigma}(p_{i})=p_{i}\circ\varsigma=\varsigma_{i}\langle p_{j}\rangle$,
therefore we can also write $\varsigma$ as $\langle p_{i}\rangle\mapsto\langle f_{\varsigma}(p_{i})\rangle$.
Conversely, for every automorphism $f:\mathbf{E}[v,0]\rightarrow\mathbf{E}[v,0]$,
define $\varsigma_f$ to be the substitution
$\langle p_{i}\rangle\mapsto\langle f(p_{i})\rangle$.
Obviously $\varsigma_f\in\mathcal{S}(v,0)$ and it follows immediately
that $\langle p_{i}\rangle\mapsto\langle f^{-1}(p_{i})\rangle$,
also in $\mathcal{S}(v,0)$, is its inverse. \bgroup \renewcommand*\theenumi{\alph{enumi}}\renewcommand*\labelenumi{\theenumi)}

Thus, \emph{$\mathcal{S}_{p}(v,0)$} is isomorphic to the group of
automorphisms of $\mathbf{E}[v,0]$. However, from lattice
theory any automorphism of a finite Boolean lattice is uniquely
determined by its values on the atoms of the lattice, and that these
automorphisms correspond to the permutations of the atoms. In our
case the atoms are the $2^{v}$ minterms of $\mathbf{E}[v,0]$, which
proves our claim.\textsc{\scriptsize{\hfill{}$\blacksquare$}}\vspace{1mm}

\begin{elabeling}{00.0000.0000.0000}
\item [{\textbf{\textsc{Theorem~4}}}] \noindent \emph{Let $\varsigma\in\mathcal{S}_{p}(v,0)$.
Then for every $\mathbf{E}[v,d]$ context:}\end{elabeling}
\begin{enumerate}
\item \emph{The function }$f_{\varsigma}:[\mathit{1}]\rightarrow[\mathit{1}]$\emph{ defined
as }$f_{\varsigma}(\mu)=\mu\circ\varsigma$\emph{ is a bijection on
the set }$_{v}^{d}[\mathit{1}]$\emph{ of minterms.}
\item \emph{The function }$f_{\varsigma}:\mathbf{E}[v,d]\rightarrow\mathbf{E}[v,d]$
\emph{defined as }$f_{\varsigma}(\varphi)=\varphi\circ\varsigma$\emph{
is a lattice automorphism.}
\end{enumerate}
\textbf{\textsc{Proof.}}\qquad{}Since $\mathbf{E}[v,d]=\wp([\mathit{1}])$
when minmatrices are viewed as sets, b) is a corollary of a).
Thus, we can prove a) by induction on the modal level $d$,
even though the induction step makes use of b).

For $d=0$ the result follows directly from Theorem 3.
Assume a) and b) hold up to some level $d$ and let $\mu$ be a level $d+1$ minterm.
By definition $\mu=\epsilon\,\psi$, where:
\begin{itemize}
\item $\epsilon$, the non-modal prefix of $\mu,$ is a level 0 minterm,
hence so is $\epsilon\circ\varsigma$.
\item $\psi$ is a product of all the level $d+1$ modal factors from the
set $\{\lozenge\phi_{k}\}$, complemented or not, with $\phi_{k}$
being all the level $d$ formulas in DCF. But by the induction hypotheses,
$\varsigma$ permutes the level $d$ minmatrices, therefore $\{\lozenge(\phi_{k}\circ\varsigma)\}=\{\lozenge\phi_{k}\}$.
Then $\psi\circ\varsigma$ is again a product of all the level $d+1$
modal factors from $\{\lozenge\phi_{k}\}$, complemented or not.
\end{itemize}

We apply property US-6 to conclude that $\mu\circ\varsigma\approx(\epsilon\circ\varsigma)(\psi\circ\varsigma)$
is a level $d+1$ minterm, so $f_{\varsigma}$ is well-defined. Also,
if $\mu_{1},\mu_{2}\in[\mathit{1}]$, then $\mu_{1}\circ\varsigma\approx\mu_{2}\circ\varsigma$
implies $\mu_{1}\circ\varsigma\circ\varsigma^{-1}\approx\mu_{2}\circ\varsigma\circ\varsigma^{-1}$
and $\mu_{1}\approx\mu_{2}$. Thus, $f_{\varsigma}$ is injective
on the finite set $[\mathit{1}]$, hence it is a bijection.\textsc{\scriptsize{\hfill{}$\blacksquare$}}\vspace{4mm}

In general, for any level 0 substitution $\sigma$ we still
have $[\mathit{1}]\circ\sigma\approx[\mathit{1}]$. Thus, the function \emph{$f_{\sigma}:[\mathit{1}]\rightarrow\mathbf{E}[v,d]$}
defined as \emph{$f_{\sigma}(\mu)=\mu\circ\sigma$ }is such that the
(possibly empty) sets $f_{\sigma}(\mu_{i})$ are disjoint (by US-6)
and $\bigcup_{\mu_{i}\in[\mathit{1}]}f_{\sigma}(\mu_{i})=[\mathit{1}]$. But $f_{\sigma}$
is a bijection on $[\mathit{1}]$ only when $\sigma$ is prime.\egroup

Theorem 4 implies that for every context $\mathbf{E}[v,d]$, $\mathcal{S}_{p}(v,0)$
determines a group action on the set $_{v}^{d}[\mathit{1}]$ of minterms.
We then define the (context-dependent) \emph{prime orbits}
of minterms as the orbits of this group action.
\vspace{1mm}

\begin{elabeling}{00.0000.0000.0000}
\item [{\textbf{\textsc{Theorem~5}}}] \emph{If a minmatrix $\varphi$ includes
some, but not all, the minterms of a prime orbit $\omega$, then it
collapses under some prime substitution $\varsigma$.}
\end{elabeling}
\textbf{\textsc{Proof.}}\qquad{}This follows from Theorem 4 and the properties of orbits.

Let $\omega^{*}=\{\mu_{1},\ldots,\mu_{k}\}$
be an incomplete, non-empty prime orbit $\omega^{*}\subsetneq\omega$ and $\mu_{x}\in\omega\setminus\omega^{*}$.
Since the restriction of a group action to an orbit is transitive,
there exists a prime substitution $\varsigma$
such that $\mu_{1}\circ\varsigma\approx\mu_{x}$. Then if $f_{\varsigma}:[\mathit{1}]\rightarrow[\mathit{1}]$
is $f_{\varsigma}(\mu)=\mu\circ\varsigma$, we have $f_{\varsigma}(\mu_{1})\notin\omega^{*}$,
so $f_{\varsigma}(\omega^{*})\not\subset\omega^{*}$. But since $\omega^{*}$
is finite, there must be a $\mu_{i}\in\omega^{*}$, $1\leq i\leq k$,
such that $\mu_{i}\notin f_{\varsigma}(\omega^{*})$. Hence for every
$\mu_{y}\in\omega^{*}$, $\mu_{y}\circ\varsigma\not\approx\mu_{i}$.

Also, given any $\mu_{y}\in\varphi\setminus\omega^{*}$ with $\mu_{y}$
in some prime orbit $\omega'\not=\omega$, we have $\mu_{y}\circ\varsigma\in\omega'$,
and because in general group action orbits are disjoint this implies $\mu_{y}\circ\varsigma\notin\omega^{*}$,
so once again $\mu_{y}\circ\varsigma\not\approx\mu_{i}$.

Thus, overall, there is no minterm $\mu_{y}\in\varphi$ such that $\mu_{y}\circ\varsigma\approx\mu_{i}$,
therefore $\mu_{i}\notin\varphi\circ\varsigma$ and $\varphi$ collapses under
the substitution $\varsigma$.\textsc{\scriptsize{\hfill{}$\blacksquare$}}

\begin{elabeling}{00.0000.0000.0000}
\item [{\textbf{\textsc{Corollary~6}}}] \noindent \emph{Every CMM must consist only of complete prime orbits}.
\end{elabeling}

So in this sense, prime orbits can be considered the ``building blocks'' of CMMs.
As mentioned before, the above is only a necessary condition
for a minmatrix $\xi$ to be the CMM of $\mathbf{S}$.
To be also a sufficient condition,
$\xi$ must not include any redundant prime orbits or minterms
(i.e.~equiprovable to $\mathit{0}$ in $\mathbf{S}$).
So $\xi=\cml\mathbf{S}\cmr$
iff $\xi$ consists precisely of all $\mathbf{S}$-satisfiable minterms.
For a concrete system $\mathbf{S}$, this is typically shown using $\mathbf{S}$-models.
But for the next results in our paper we do not need to use this method,
as we only rely on the fact that
every ${}_{v}^{d}\cml\mathbf{S}\cmr$ exists,
due to the finiteness of the context.

For an $\mathbf{E}[v,d]$ system $\mathbf{S}$, or for its CMM $\cml\mathrm{\mathbf{S}}\cmr$,
let $\Omega(\mathbf{S})=\Omega(_{v}^{d}\mathbf{S})$
be the set of prime orbits from $_{v}^{d}\cml\mathbf{S}\cmr$. Obviously,
$\Omega(\mathbf{E})=\Omega([\mathit{1}])$ in all contexts.

\pagebreak{}

\paragraph{Examples.}

Figure 1 shows the lattice of $\mathbf{E}[0,1]$ systems and their
CMMs. It has $\mathbf{E}$ as top and $\mathbf{F}$ as bottom. The
other systems are named for the purpose of this example only, as
they do not represent the normal systems usually denoted by the letters.
(However, $\mathbf{K}$ and $\mathbf{D}$ do have some connection with
their normal counterparts, as explained below.)\vspace{-2mm}

\begin{figure}[H]
\begin{centering}
\includegraphics[scale=0.63]{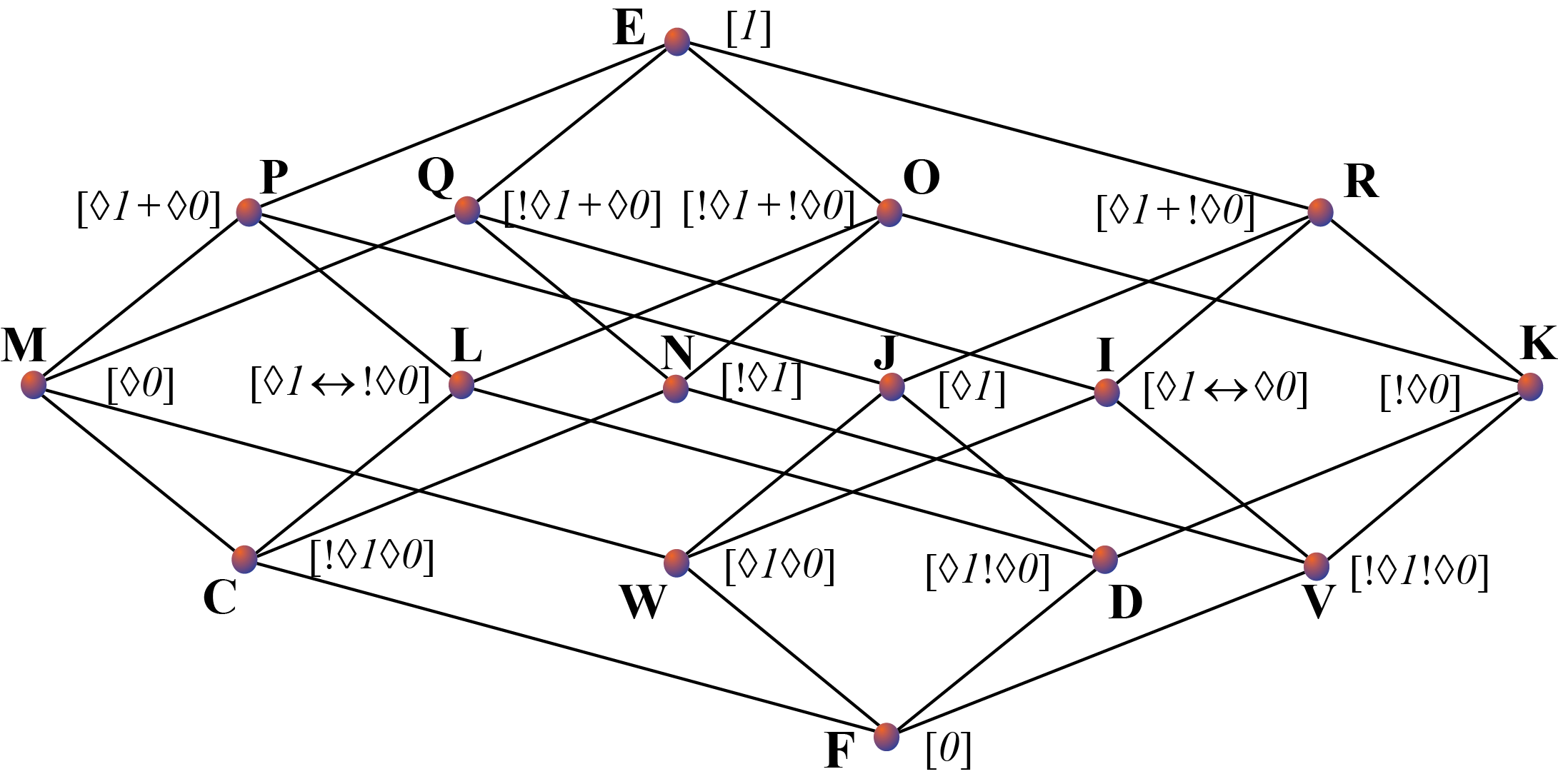}
\par\end{centering}

\caption{\label{fig:E01-systems} The $\mathbf{E}[0,1]$ systems and
their CMMs}
\end{figure}\vspace{-2mm}

Considering axiom $\mathrm{D}=\lozenge\mathit{1}\mbox{!}\lozenge\mathit{0}$,
this minmatrix has a single minterm, so it could only collapse
to $[\mathit{0}]$. But this is not the case, otherwise the normal
modal system $\mathbf{D}$ would collapse to $\mathbf{F}$. The other
atoms of the above lattice will be shown to be non-collapsing CMMs
too by Theorem 14. The corresponding axioms are disjoint formulas,
so these atoms yield 4 distinct $\mathbf{E}[0,1]$ systems. Next,
we use Theorem 2 a) to see that all the other elements in the lattice
are indeed CMMs, so there are 16 distinct $\mathbf{E}[0,1]$ systems.

In subsequent contexts this lattice is ``fractally refined''
by new systems that are not $\mathbf{E}[0,1]$ systems.
For example, the normal system $\mathbf{Ver}$
is obtained by adding to $\mathbf{E}$
the axiom $\mbox{!}\lozenge p$,
which makes it an $\mathbf{E}[1,1]$ system.
Since $\mathrm{V}=\mbox{!}\lozenge\mathit{1}\mbox{!}\lozenge\mathit{0}$
is a theorem of $\mathbf{Ver}$
we have $\cml\mathbf{Ver}\cmr\subset\cml\mathbf{V}\cmr$.
Then the normal systems $\mathbf{K}$ and $\mathbf{D}$, when defined
as extensions of $\mathbf{E}$, are $\mathbf{E}[2,1]$ systems (because
of axiom K).
Their $\mathbf{E}[2,1]$ CMMs are similarly shown to be
included in those of systems $\mathbf{K}$ and $\mathbf{D}$ from the above
diagram respectively.

In fact, using the models presented in \cite{Fine}, one can show
that the $\mathbf{E}[2,1]$ CMM of the normal system $\mathbf{K}$
has 64 minterms, where the state of $\lozenge\mathit{0}$ is always
0, then $p$, $q$, $\lozenge(pq)$, $\lozenge(p!q)$, $\lozenge(!pq)$ and $\lozenge(!p!q)$
have all 64 states combinations, while the states of the remaining
$\mathbf{E}[2,1]$ modal factors are uniquely determined from
$\lozenge(\varphi+\psi)\approx\lozenge\varphi+\lozenge\psi$
(e.g.~$\lozenge p\approx\lozenge(pq)+\lozenge(p!q)$).
This also shows how DNF minterms are a special case of DCF minterms.

It is also instructive to present here an example of a CMM intersection
collapse for Theorem 2 b), as this situation will be mentioned in Theorems
14 and 15. A quite relevant example can be given for base $\mathbf{K}$,
where our theory of CMMs and prime orbits also applies, as we shall
see from \cite{Soncodi}, with the modal DNF instead of the DCF. However,
since determining CMMs in concrete cases is quite laborious, we can
only sketch the proof here.

The example involves some well-known extensions of the normal system
$\mathbf{T}$, namely $\mathbf{S4}$, $\mathbf{B}$ and $\mathbf{S5}$.
The 9 DNF factors in $\mathbf{K}[1,2]$ are
$p$, $\lozenge(p\lozenge p\lozenge!p)$, $\lozenge(p\lozenge p!\lozenge!p)$,
$\lozenge(p!\lozenge p\lozenge!p)$, $\lozenge(p!\lozenge p!\lozenge!p)$,
$\lozenge(!p\lozenge p\lozenge!p)$, $\lozenge(!p\lozenge p!\lozenge!p)$, 
$\lozenge(!p!\lozenge p\lozenge!p)$, $\lozenge(!p!\lozenge p!\lozenge!p)$.
But as the $\mathbf{K}[1,1]$ CMM of $\mathbf{T}$
can be shown to be $\varphi$ from (\ref{eq:CMM-T}),
the modal factors
$\lozenge(p!\lozenge p\lozenge!p)$, $\lozenge(p!\lozenge p!\lozenge!p)$,
$\lozenge(!p\lozenge p!\lozenge!p)$, $\lozenge(!p!\lozenge p!\lozenge!p)$
are equiprovable to $\mathit{0}$ in all extensions of $\mathbf{T}$,
so their state must be 0 in all the minterms of those CMMs.
Consequently we can omit the corresponding rows from the
minmatrix representation without any loss of information.
Then using the $\mathbf{S4}$ model from \cite{Moss} one can show that
in context $\mathbf{K}[1,2]$:\vspace{-3.5mm}

\[
_{1}^{2}\cml\mathbf{S4}\cmr=\begin{array}{c|cccccccccccccc|}
p & 1 & 1 & 1 & 1 & 1 & 1 & 1 & 0 & 0 & 0 & 0 & 0 & 0 & 0\\
\hline \lozenge(\phantom{!}p\phantom{!}\lozenge p\phantom{!}\lozenge!p) & 1 & 1 & 1 & 1 & 1 & 1 & 0 & 1 & 1 & 1 & 1 & 0 & 0 & 0\\
\lozenge(\phantom{!}p\phantom{!}\lozenge p!\lozenge!p) & 1 & 1 & 1 & 0 & 0 & 0 & 1 & 1 & 1 & 0 & 0 & 1 & 1 & 0\\
\lozenge(!p\phantom{!}\lozenge p\phantom{!}\lozenge!p) & 1 & 1 & 0 & 1 & 1 & 0 & 0 & 1 & 1 & 1 & 1 & 1 & 1 & 0\\
\lozenge(!p!\lozenge p\phantom{!}\lozenge!p) & 1 & 0 & 1 & 1 & 0 & 1 & 0 & 1 & 0 & 1 & 0 & 1 & 0 & 1
\end{array}
\]\vspace{0.3mm}
and from a similar model for $\mathbf{B}$ we can obtain:\vspace{-2mm}

\[
_{1}^{2}\cml\mathbf{B}\cmr=\begin{array}{c|cccccccc|}
p & 1 & 1 & 1 & 1 & 0 & 0 & 0 & 0\\
\hline \lozenge(\phantom{!}p\phantom{!}\lozenge p\phantom{!}\lozenge!p) & 1 & 1 & 1 & 0 & 1 & 1 & 0 & 0\\
\lozenge(\phantom{!}p\phantom{!}\lozenge p!\lozenge!p) & 1 & 1 & 0 & 1 & 0 & 0 & 0 & 0\\
\lozenge(!p\phantom{!}\lozenge p\phantom{!}\lozenge!p) & 1 & 0 & 1 & 0 & 1 & 1 & 1 & 0\\
\lozenge(!p!\lozenge p\phantom{!}\lozenge!p) & 1 & 0 & 0 & 0 & 1 & 0 & 1 & 1
\end{array}
\]\vspace{0.1mm}

The intersection $[\mathrm{B4}]={}_{1}^{2}\cml\mathbf{S4}\cmr\cap{}_{1}^{2}\cml\mathbf{B}\cmr$
has 6 minterms.
But this is \emph{not} $_{1}^{2}\cml\mathbf{S5}\cmr$, since it can be collapsed,
using the $\mathbf{S5}$ theorem $\lozenge\lozenge p \rightarrow \square\lozenge p$,
to:\vspace{-1.8mm}

\[
_{1}^{2}\cml\mathbf{S5}\cmr=\begin{array}{c|cccc|}
p & 1 & 1 & 0 & 0\\
\hline \lozenge(\phantom{!}p\phantom{!}\lozenge p\phantom{!}\lozenge!p) & 1 & 0 & 1 & 0\\
\lozenge(\phantom{!}p\phantom{!}\lozenge p!\lozenge!p) & 0 & 1 & 0 & 0\\
\lozenge(!p\phantom{!}\lozenge p\phantom{!}\lozenge!p) & 1 & 0 & 1 & 0\\
\lozenge(!p!\lozenge p\phantom{!}\lozenge!p) & 0 & 0 & 0 & 1
\end{array}
\]\vspace{0.4mm}

The latter no longer collapses,
as the only CMM strictly included in it is $_{1}^{2}\cml\mathbf{Triv}\cmr$.
Thus, even though $[\mathrm{B4}]$ consists only of complete prime
orbits (here, pairs of minterms), it cannot be the CMM of any system,
since as a theorem it alone collapses any such system to $\mathbf{S5}$.

\pagebreak{}

\section{Automorphisms of $\mathbf{E}\cml v,d\cmr$}

\paragraph*{Uniform replacements.}

In this paper we define a set of formula transformations called uniform
replacements (UR). These are similar to uniform substitutions, only
instead of replacing all the occurrences of propositional variables
in a formula we replace all the occurrences of certain sub-formulas.

We shall focus on a particular class of context-preserving uniform
replacements associated with formulas of modal degree 1 in one variable,
$\rho=\rho(e)$. Such a formula can also be written as $\rho(e)=\eta(e,\Diamond e,\lozenge!e)$,
or alternately $\rho(e)=\eta(e,\Diamond e,\square e)$, for a suitable
level 0 formula $\eta=\eta(a,b,c)$. Up to algebraic equivalence,
there are 256 level 0 formulas in 3 variables and thus 256 formulas
$\rho$, for example $\rho(e)=e\lozenge e+\mbox{!}e\square e$ or
$\rho(e)=e\leftrightarrow\lozenge e$.

Let $\rho=\rho(e)$ be such a level 1 formula. The result, denoted
as $\varphi*\rho$, of applying the \emph{uniform replacement} corresponding
to $\rho$ to a formula $\varphi=\varphi\langle p_{i}\rangle$ is
the formula $(\varphi*\rho)\langle p_{i}\rangle$ obtained by recursively
replacing the sub-formulas $\psi$, $\theta$, $\ldots$ that occur
in $\varphi$ according to the following rules:
\begin{elabeling}{00.00.0000}
\item [{(UR-1)}] $\mathit{0}*\rho\triangleq\mathit{0}$ and $\mathit{1}*\rho\triangleq\mathit{1}$
\item [{(UR-2)}] $p_{i}*\rho\triangleq p_{i}$
\item [{(UR-3)}] $(\psi+\theta)*\rho\triangleq\psi*\rho+\theta*\rho$
\item [{(UR-4)}] $(!\psi)*\rho\triangleq\mbox{!}(\psi*\rho)$
\item [{(UR-5)}] $(\Diamond\psi)*\rho\triangleq\rho(\psi*\rho)$
\end{elabeling}

Observe that we equate the UR with its defining level 1 formula $\rho$,
so the distinction must to be made from the context. Some immediate
consequences of the above definitions are:
\begin{elabeling}{00.00.0000}
\item [{(UR-6)}] $(\psi\,\theta)*\rho\approx(\psi*\rho)(\theta*\rho)$
\item [{(UR-7)}] $(\psi\rightarrow\theta)*\rho\approx(\psi*\rho\rightarrow\theta*\rho)$
\item [{(UR-8)}] $(\psi\leftrightarrow\theta)*\rho\approx(\psi*\rho\leftrightarrow\theta*\rho)$
\item [{(UR-9)}] $(\psi\not\leftrightarrow\theta)*\rho\approx(\psi*\rho\not\leftrightarrow\theta*\rho)$
\item [{(UR-10)}] $(\square\psi)*\rho\approx\mbox{!}\rho(!(\psi*\rho))$
\end{elabeling}

Note the recursive application of $\rho$ to the sub-formulas that
involve the modal operators.

As defined, a UR transformation $\rho$ is \emph{context-preserving},
i.e.~it does not increase the level or number of variables in the
formula. This transformation leaves a level 0 formula $e$ unchanged.
And for a level 1 formula $\varphi$, it consists of the following
uniform replacements in $\varphi$: 
\begin{elabeling}{00.00.0000}
\item [{(UR-$\Diamond$)}] $\Diamond e\:\mapsto\:\rho(e)$
\item [{(UR-$\square$)}] $\square e\:\mapsto\:\mbox{!}\rho(!e)$
\end{elabeling}
where the distinct level 0 sub-formulas $e$ that actually occur
under the modal operators in $\varphi$ result in different replacements
in $\varphi$, but uniform in case of multiple occurrences of the
same $e$ under several instances of the modal operators. In this
case UR-$\lozenge$ and UR-$\square$ can be combined into:
\begin{elabeling}{00.00.0000}
\item [{(UR-$\rho$)}] $(\Diamond e,\square e)\:\mapsto\:(\,\rho(e),\,\mbox{!}\rho(!e)\,)$
\end{elabeling}

Alternately, we could require the textual replacement of $\mbox{`}\square\mbox{\mbox{'}}$
with $\mbox{`}!\lozenge!\mbox{'}$ in $\varphi$ prior to applying
$\rho$, in which case UR-$\square$ becomes redundant.

Also note that UR-$\rho$ applies even to constant modalities, if
any, hence $\lozenge\mathit{1}\:\mapsto\:\rho(\mathit{1})$ and $\lozenge\square\mathit{0}\:\mapsto\:\rho((\square\mathit{0})*\rho)\approx\rho(!\rho(!\mathit{0}*\rho))\approx\rho(!\rho(\mathit{1}))$
etc.  \bgroup \renewcommand*\theenumi{\alph{enumi}}\renewcommand*\labelenumi{\theenumi)}\vspace{2mm}

\begin{elabeling}{00.0000.0000.0000}
\item [{\textbf{\textsc{Lemma~7}}}] \emph{Let }$\varphi$\emph{ and }$\psi$\emph{
be any formulas. Then: }\end{elabeling}
\begin{enumerate}
\item \emph{If }$\vdash\mathbf{_{\hspace{-3pt}_{E}}\,}\varphi$\emph{ then
}$\vdash\mathbf{_{\hspace{-3pt}_{E}}\,}\varphi*\rho$.
\item \emph{If }$\varphi\approx\psi$\emph{ then }$\varphi*\rho\approx\psi*\rho$.
\item \emph{If }$\varphi\,\psi\approx\mathit{0}$\emph{ then }$(\varphi*\rho)(\psi*\rho)\approx\mathit{0}$.
\end{enumerate}
\textbf{\textsc{Proof.}}\qquad{}From the UR properties we have $\varphi*\rho\approx[\varphi]*\rho$.
But $\vdash\mathbf{_{\hspace{-3pt}_{E}}\,}\varphi$ implies $[\varphi]=e(\langle p_{i}\rangle,\langle\lozenge\phi_{j}\langle p_{k}\rangle\rangle)$
for some level 0 formula $e(\langle p_{i}\rangle,\langle q_{j}\rangle)\approx\mathit{1}$,
i.e.~when the modal factors $\lozenge\phi_{j}\langle p_{k}\rangle$
are considered atoms, $e$ is a tautology. Then $\varphi*\rho\approx[\varphi]*\rho\approx e(\langle p_{i}\rangle,\langle\lozenge\phi_{j}\langle p_{k}\rangle*\rho\rangle)\approx e(\langle p_{i}\rangle,\langle r_{j}\rangle)\approx\mathit{1}$,
therefore $\vdash\mathbf{_{\hspace{-2pt}_{E}}\,}\varphi*\rho$, which
proves a). Since $\varphi\approx\psi$ is $\,\vdash\mathbf{_{\hspace{-3pt}_{E}}\,}\varphi\leftrightarrow\psi$,
b) follows from a) and UR-8. Lastly, c) follows from b), UR-1 and
UR-6.\textsc{\scriptsize{\hfill{}$\blacksquare$}}\vspace{3mm}

Note that a) holds in $\mathbf{E}$, but there is no ``rule of UR''
for an arbitrary system $\mathbf{S}$: if $\vdash\mathbf{_{\hspace{-3pt}_{\mathbf{S}}}\,}\varphi$,
then typically $\not\vdash\mathbf{_{\hspace{-3pt}_{\mathbf{S}}}\,}\varphi*\rho$.
(In fact, we will show that $\varphi*\rho$ is a theorem of the transformed
system $\mathbf{S}_{\rho}\triangleq\mathbf{S}*\rho$.)

Given a sub-formula $\psi$, we write $\text{¡}\psi$ when either
$\psi$ or $!\psi$ occurs at some position in a formula $\varphi$.
Then the generic DCF of any level $d$ formula $\varphi$ is:

\vspace{-1mm}

\begin{equation}
\varphi\:\approx\:\sum_{k}\epsilon_{k}\prod_{l}\text{¡}\Diamond\phi_{l}\:=\:\sum_{k}\epsilon_{k}\langle p_{i}\rangle\prod_{l}\text{¡}\Diamond\phi_{l}\langle p_{i}\rangle\label{eq:DCFform}
\end{equation}
where $k$ and $l$ have suitable ranges, $\epsilon_{k}=\prod_{i}\text{¡}p_{i}$
are the minterm prefixes (the level 0 minterms; for $v=0$ there is
a single minterm, namely~$\epsilon=\mathit{1}$) and $\Diamond\phi_{l}$ are the
level $d$ modal factors, which range over all $\phi_{l}\in\mathbf{E}[v,d-1]$.\egroup

If $\sigma=\langle\sigma_{i}\rangle$ is a substitution, let $\sigma*\rho\triangleq\langle\sigma_{i}*\rho\rangle$.
Then for any prefix $\epsilon$:
\begin{equation}
(\epsilon\circ\sigma)*\rho\approx(\prod_{i}\text{¡}\sigma_{i})*\rho\approx\prod_{i}\text{¡}(\sigma_{i}*\rho)\approx\epsilon\circ(\sigma*\rho)\approx(\epsilon*\rho)\circ(\sigma*\rho)\label{eq:UR-sigma}
\end{equation}

\pagebreak{}

\begin{elabeling}{00.0000.0000.0000}
\item [{\textbf{\textsc{Lemma~8}}}] \emph{Let $\varphi$ be a formula,
}$\sigma$\emph{ a substitution and $\rho$ a UR. Then:}\vspace{-1mm}
\end{elabeling}
\[
(\varphi\circ\sigma)*\rho\approx(\varphi*\rho)\circ(\sigma*\rho)
\]
\textbf{\textsc{Proof.}}\qquad{}We proceed by induction on the modal
degree $d$ of $\varphi$. For $d=0$, $\varphi=\sum_{k}\epsilon_{k}$
and the result follows immediately from (\ref{eq:UR-sigma}). Assuming
the claim is true up to level $d$, we apply the defining US and UR
properties to a level $d+1$ formula $\varphi$ expanded as in (\ref{eq:DCFform})
to obtain: \vspace{-5mm}

\begin{eqnarray*}
(\varphi\circ\sigma)*\rho & \approx & \left(\sum_{k}(\epsilon_{k}\circ\sigma)\prod_{l}\text{¡}\lozenge(\phi_{l}\circ\sigma)\right)*\rho\\
 & \approx & \sum_{k}((\epsilon_{k}\circ\sigma)*\rho)\prod_{l}(\text{¡}\lozenge(\phi_{l}\circ\sigma)*\rho)\\
 & \approx & \sum_{k}((\epsilon_{k}*\rho)\circ(\sigma*\rho))\prod_{l}\text{¡}\rho((\phi_{l}\circ\sigma)*\rho)\\
 & \approx & \sum_{k}((\epsilon_{k}*\rho)\circ(\sigma*\rho))\prod_{l}\text{¡}\rho((\phi_{l}*\rho)\circ(\sigma*\rho))\\
 & \approx & \sum_{k}((\epsilon_{k}*\rho)\circ(\sigma*\rho))\prod_{l}(\text{¡}\rho(\phi_{l}*\rho)\circ(\sigma*\rho))\\
 & \approx & \left(\sum_{k}(\epsilon_{k}*\rho)\prod_{l}(\text{¡}\lozenge\phi_{l}*\rho)\right)\circ(\sigma*\rho)\\
 & \approx &
 (\varphi*\rho)\circ(\sigma*\rho)
\end{eqnarray*}
because $\phi_{l}$ are formulas of level $\leq d$. The above makes use of
(\ref{eq:UR-sigma}) on lines 2-3 and of $\rho(\varphi\circ\sigma)\approx\rho(\varphi)\circ\sigma$
on lines 4-5.  The latter is obvious when we observe that $\rho(\varphi)\approx\rho\circ\varphi$ as a 1-variable substitution in $\rho(e)$.\textsc{\scriptsize{\hfill{}$\blacksquare$}}\vspace{4mm}

Let $\mathcal{R}(1)$ be the set of all URs up to DCF equivalence,
so $\mathcal{R}(1)\subset\mathbf{E}[1,1]$, with $|\mathcal{R}(1)|=256$
distinct minmatrices. Observe that URs were defined independently
of the formulas $\varphi$ from other $\mathbf{E}[v,d]$ contexts.
The reason is that we intend to define a monoid structure on $\mathcal{R}(1)$,
and then also a (right) monoid action $A:\mathbf{E}[v,d]\times\mathcal{R}(1)\rightarrow\mathbf{E}[v,d]$,
with $A(\varphi,\rho)\triangleq\varphi*\rho$.

To achieve this, we define the \emph{composition} $\rho\rho'$ of
URs $\rho$ and $\rho'$ as:\vspace{-2mm}

\[
\rho\rho'\triangleq\rho*\rho'
\]
where the right hand side $\rho$ is taken to be the associated level
1 formula $\rho(e)=\eta(e,\Diamond e,\lozenge!e)\in\mathbf{E}[1,1]$,
such that $\rho'$ can be applied to it.

We now verify that the above are well-defined. The fact that $\varrho_{0}(e)=\Diamond e$
(i.e.~$\Diamond e\mapsto\Diamond e$) is the monoid unit in $\mathcal{R}(1)$
is straightforward to check. For the monoid action, the identity axiom
is proven by the result below.\vspace{3mm}

\begin{elabeling}{00.0000.0000.0000}
\item [{\textbf{\textsc{Lemma~9}}}] \emph{For any formula} $\varphi$,\emph{
we have }$\varphi*\varrho_{0}\approx\varphi$.
\end{elabeling}
\noindent \textbf{\textsc{Proof.}}\qquad{}We write $\varphi\langle p_{i}\rangle=\delta(\langle p_{i}\rangle,\langle\lozenge\phi_{j}\langle p_{k}\rangle\rangle)$,
with $\delta$ a level 0 formula, which yields:\vspace{2mm}

\noindent 
\[
\varphi*\varrho_{0}\:\approx\:\delta(\langle p_{i}\rangle,\langle\lozenge\phi_{j}\rangle)*\varrho_{0}\:\approx\:\delta(\langle p_{i}\rangle,\langle\varrho_{0}(\phi_{j}*\varrho_{0})\rangle)\:\approx\:\delta(\langle p_{i}\rangle,\langle\lozenge(\phi_{j}*\varrho_{0})\rangle)
\]
and the result follows immediately by induction on the level
$d$ of $\varphi$.\hfill{}\textsc{\scriptsize{$\blacksquare$}}{\scriptsize \par}

\vspace{5mm}
The next lemma can be used to establish both the associativity of
UR composition on $\mathcal{R}$(1) and the compatibility condition
required by $A$ for it:\vspace{3mm}

\begin{elabeling}{00.0000.0000.0000}
\item [{\textbf{\textsc{Lemma~10}}}] \emph{Let }$\varphi$\emph{ be any
formula and $\rho,\rho'\in\mathcal{R}(1)$. Then}:\vspace{-4mm}
\end{elabeling}
\begin{eqnarray*}
(\varphi*\rho)*\rho' & \approx & \varphi*(\rho\rho')
\end{eqnarray*}
\noindent \textbf{\textsc{Proof.}}\qquad{}The proof is by induction
on the level $d$ of $\varphi$. For $d=0$, the result is obvious,
since $\varphi$ is unchanged by URs. Assume the claim holds up to
some level $d$ and let $\varphi$ be a level $d+1$ formula. Let
$\delta$ be the level 0 formula such that $\varphi\langle p_{i}\rangle=\delta(\langle p_{i}\rangle,\langle\lozenge\phi_{j}\langle p_{k}\rangle\rangle)$,
where all $\phi_{j}$ are of level $\leq d$. Then:
\setlength{\jot}{2mm}\vspace{-5mm}

\begin{eqnarray*}
(\varphi*\varrho)*\varrho' & \approx & \delta(\langle p_{i}\rangle,\langle\lozenge\phi_{j}*\rho\rangle)*\rho'\\
 & \approx & \delta(\langle p_{i}\rangle,\langle\rho(\phi_{j}*\rho)\rangle)*\rho'\\
 & \approx & \delta(\langle p_{i}\rangle,\langle\eta(\phi_{j}*\rho,\lozenge(\phi_{j}*\rho),\lozenge!(\phi_{j}*\rho))\rangle)*\rho'\\
 & \approx & \delta(\langle p_{i}\rangle,\langle\eta((\phi_{j}*\rho)*\rho',\rho'((\phi_{j}*\rho)*\rho'),\rho'(!(\phi_{j}*\rho)*\rho')))\rangle)\\
 & \approx & \delta(\langle p_{i}\rangle,\langle\eta(\alpha_{j},\rho'(\alpha_{j}),\rho'(!\alpha_{j}))\rangle)\\
\varphi*(\rho\rho') & \approx & \delta(\langle p_{i}\rangle,\langle\lozenge\phi_{j}\rangle)*(\rho\rho')\\
 & \approx & \delta(\langle p_{i}\rangle,\langle\lozenge\phi_{j}\rangle)*\eta(e,\rho'(e),\rho'(!e))\\
 & \approx & \delta(\langle p_{i}\rangle,\langle\eta(\phi_{j}*(\rho\rho'),\rho'(\phi_{j}*(\rho\rho')),\rho'(!\phi_{j}*(\rho\rho')))\rangle)\\
 & \approx & \delta(\langle p_{i}\rangle,\langle\eta(\beta_{j},\rho'(\beta_{j}),\rho'(!\beta_{j}))\rangle)
\end{eqnarray*}
where $\alpha_{j}=(\phi_{j}*\rho)*\rho'$ and $\beta_{j}=\phi_{j}*(\rho\rho')$.
But since $\phi_{j}$ are of level $\leq d$ and the URs are context-preserving,
$\alpha$ and $\beta$ must also be of level $\leq d$, hence by the
induction step $\alpha\approx\beta$, which proves our claim.\textsc{\scriptsize{\hfill{}$\blacksquare$}}\vspace{4mm}

\noindent \textbf{Prime URs.} Among the 256 URs there are 24 that
are invertible, and we call them \emph{prime UR transformations}.

Let $\mathcal{R}_{p}(1)=\{\varrho_{i}:0\leq i<24\}\subset\mathcal{R}(1)$
be the set of prime URs, which are defined in Table 2.
These correspond to the 24 permutations
from the symmetric group $\mathrm{S}_4$ as per column 4,
which we prove next.

\pagebreak{}

\noindent \begin{center}
\begin{tabular}{|c|l|l|l|c|}
\multicolumn{5}{c}{{\footnotesize{Table 2: List of prime UR transformations}}}\tabularnewline[2mm]
\hline 
\multicolumn{1}{|c|}{{\scriptsize{Prime}}} & {\scriptsize{$\qquad\varrho_{i}(e)$ }} & {\scriptsize{$\qquad!\varrho_{i}(!e)$ }} & {\scriptsize{$\mathrm{S_{4}}$ correspondent}} & {\scriptsize{Inverse}}\tabularnewline
{\scriptsize{UR $\varrho_{i}$}} & {\scriptsize{\quad{}i.e. $\Diamond e\mapsto$}} & {\scriptsize{\quad{}i.e. $\square e\mapsto$}} & {\scriptsize{$(\mathrm{W,D,C,V})\mapsto$}} & {\scriptsize{UR $\varrho_{i}^{-1}$}}\tabularnewline
\hline 
\noalign{\vskip1pt}
{\scriptsize{$\varrho_{0}$}} & {\scriptsize{$\phantom{!}\Diamond\phantom{!}e$}} & {\scriptsize{$\phantom{!}\square\phantom{!}e$}} & {\scriptsize{$(\mathrm{W,D,C,V})$}} & {\scriptsize{$\varrho_{0}$}}\tabularnewline
{\scriptsize{$\varrho_{1}$}} & {\scriptsize{$\phantom{!}\Diamond\phantom{!}e\leftrightarrow\phantom{!}e$}} & {\scriptsize{$\phantom{!}\square\phantom{!}e\leftrightarrow\mbox{!}e$}} & {\scriptsize{$(\mathrm{D,W,V,C})$}} & {\scriptsize{$\varrho_{1}$}}\tabularnewline
{\scriptsize{$\varrho_{2}$}} & {\scriptsize{$\phantom{!}\Diamond\phantom{!}e\leftrightarrow\phantom{!}e+\phantom{!}\Diamond!e$}} & {\scriptsize{$\phantom{!}\square\phantom{!}e\leftrightarrow\mbox{!}e+\mbox{!}\square!e$}} & {\scriptsize{$(\mathrm{W,D,V,C})$}} & {\scriptsize{$\varrho_{2}$}}\tabularnewline
{\scriptsize{$\varrho_{3}$}} & {\scriptsize{$\phantom{!}\Diamond\phantom{!}e\leftrightarrow\phantom{!}e+\mbox{!}\Diamond!e$}} & {\scriptsize{$\phantom{!}\square\phantom{!}e\leftrightarrow\mbox{!}e+\phantom{!}\square!e$}} & {\scriptsize{$(\mathrm{D,W,C,V})$}} & {\scriptsize{$\varrho_{3}$}}\tabularnewline
{\scriptsize{$\varrho_{4}$}} & {\scriptsize{$\phantom{!}\Diamond\phantom{!}e\leftrightarrow\mbox{!}e+\phantom{!}\Diamond!e$}} & {\scriptsize{$\phantom{!}\square\phantom{!}e\leftrightarrow\phantom{!}e+\mbox{!}\square!e$}} & {\scriptsize{$(\mathrm{W,V,C,D})$}} & {\scriptsize{$\varrho_{4}$}}\tabularnewline
{\scriptsize{$\varrho_{5}$}} & {\scriptsize{$\phantom{!}\Diamond\phantom{!}e\leftrightarrow\mbox{!}e+\mbox{!}\Diamond!e$}} & {\scriptsize{$\phantom{!}\square\phantom{!}e\leftrightarrow\phantom{!}e+\phantom{!}\square!e$}} & {\scriptsize{$(\mathrm{C,D,W,V})$}} & {\scriptsize{$\varrho_{5}$}}\tabularnewline
\hline 
\noalign{\vskip1pt}
{\scriptsize{$\varrho_{6}$}} & {\scriptsize{$!\Diamond\phantom{!}e$}} & {\scriptsize{$!\square\phantom{!}e$}} & {\scriptsize{$(\mathrm{V,C,D,W})$}} & {\scriptsize{$\varrho_{6}$}}\tabularnewline
{\scriptsize{$\varrho_{7}$}} & {\scriptsize{$!\Diamond\phantom{!}e\leftrightarrow\phantom{!}e$}} & {\scriptsize{$!\square\phantom{!}e\leftrightarrow\mbox{!}e$}} & {\scriptsize{$(\mathrm{C,V,W,D})$}} & {\scriptsize{$\varrho_{7}$}}\tabularnewline
{\scriptsize{$\varrho_{8}$}} & {\scriptsize{$!\Diamond\phantom{!}e\leftrightarrow\phantom{!}e+\phantom{!}\Diamond!e$}} & {\scriptsize{$!\square\phantom{!}e\leftrightarrow\mbox{!}e+\mbox{!}\square!e$}} & {\scriptsize{$(\mathrm{C,V,D,W})$}} & {\scriptsize{$\varrho_{9}$}}\tabularnewline
{\scriptsize{$\varrho_{9}$}} & {\scriptsize{$!\Diamond\phantom{!}e\leftrightarrow\phantom{!}e+\mbox{!}\Diamond!e$}} & {\scriptsize{$!\square\phantom{!}e\leftrightarrow\mbox{!}e+\phantom{!}\square!e$}} & {\scriptsize{$(\mathrm{V,C,W,D})$}} & {\scriptsize{$\varrho_{8}$}}\tabularnewline
{\scriptsize{$\varrho_{10}$}} & {\scriptsize{$!\Diamond\phantom{!}e\leftrightarrow\mbox{!}e+\phantom{!}\Diamond!e$}} & {\scriptsize{$!\square\phantom{!}e\leftrightarrow\phantom{!}e+\mbox{!}\square!e$}} & {\scriptsize{$(\mathrm{D,C,V,W})$}} & {\scriptsize{$\varrho_{11}$}}\tabularnewline
{\scriptsize{$\varrho_{11}$}} & {\scriptsize{$!\Diamond\phantom{!}e\leftrightarrow\mbox{!}e+\mbox{!}\Diamond!e$}} & {\scriptsize{$!\square\phantom{!}e\leftrightarrow\phantom{!}e+\phantom{!}\square!e$}} & {\scriptsize{$(\mathrm{V,W,D,C})$}} & {\scriptsize{$\varrho_{10}$}}\tabularnewline
\hline 
\noalign{\vskip1pt}
{\scriptsize{$\varrho_{12}$}} & {\scriptsize{$\phantom{!}\Diamond!e$}} & {\scriptsize{$\phantom{!}\square!e$}} & {\scriptsize{$(\mathrm{W,C,D,V})$}} & {\scriptsize{$\varrho_{12}$}}\tabularnewline
{\scriptsize{$\varrho_{13}$}} & {\scriptsize{$\phantom{!}\Diamond!e\leftrightarrow\phantom{!}e$}} & {\scriptsize{$\phantom{!}\square!e\leftrightarrow\mbox{!}e$}} & {\scriptsize{$(\mathrm{C,W,V,D})$}} & {\scriptsize{$\varrho_{19}$}}\tabularnewline
{\scriptsize{$\varrho_{14}$}} & {\scriptsize{$\phantom{!}\Diamond!e\leftrightarrow\phantom{!}e+\phantom{!}\Diamond e$}} & {\scriptsize{$\phantom{!}\square!e\leftrightarrow\mbox{!}e+\mbox{!}\square e$}} & {\scriptsize{$(\mathrm{W,C,V,D})$}} & {\scriptsize{$\varrho_{16}$}}\tabularnewline
{\scriptsize{$\varrho_{15}$}} & {\scriptsize{$\phantom{!}\Diamond!e\leftrightarrow\phantom{!}e+\mbox{\ensuremath{!}}\Diamond e$}} & {\scriptsize{$\phantom{!}\square!e\leftrightarrow\mbox{!}e+\phantom{!}\square e$}} & {\scriptsize{$(\mathrm{C,W,D,V})$}} & {\scriptsize{$\varrho_{17}$}}\tabularnewline
{\scriptsize{$\varrho_{16}$}} & {\scriptsize{$\phantom{!}\Diamond!e\leftrightarrow\mbox{!}e+\phantom{!}\Diamond e$}} & {\scriptsize{$\phantom{!}\square!e\leftrightarrow\phantom{!}e+\mbox{!}\square e$}} & {\scriptsize{$(\mathrm{W,V,D,C})$}} & {\scriptsize{$\varrho_{14}$}}\tabularnewline
{\scriptsize{$\varrho_{17}$}} & {\scriptsize{$\phantom{!}\Diamond!e\leftrightarrow\mbox{!}e+\mbox{\mbox{!}}\Diamond e$}} & {\scriptsize{$\phantom{!}\square!e\leftrightarrow\phantom{!}e+\phantom{!}\square e$}} & {\scriptsize{$(\mathrm{D,C,W,V})$}} & {\scriptsize{$\varrho_{15}$}}\tabularnewline
\hline 
\noalign{\vskip1pt}
{\scriptsize{$\varrho_{18}$}} & {\scriptsize{$!\Diamond!e$}} & {\scriptsize{$!\square!e$}} & {\scriptsize{$(\mathrm{V,D,C,W})$}} & {\scriptsize{$\varrho_{18}$}}\tabularnewline
{\scriptsize{$\varrho_{19}$}} & {\scriptsize{$!\Diamond!e\leftrightarrow\phantom{!}e$}} & {\scriptsize{$!\square!e\leftrightarrow\mbox{!}e$}} & {\scriptsize{$(\mathrm{D,V,W,C})$}} & {\scriptsize{$\varrho_{13}$}}\tabularnewline
{\scriptsize{$\varrho_{20}$}} & {\scriptsize{$!\Diamond!e\leftrightarrow\phantom{!}e+\phantom{!}\Diamond e$}} & {\scriptsize{$!\square!e\leftrightarrow\mbox{!}e+\mbox{!}\square e$}} & {\scriptsize{$(\mathrm{D,V,C,W})$}} & {\scriptsize{$\varrho_{23}$}}\tabularnewline
{\scriptsize{$\varrho_{21}$}} & {\scriptsize{$!\Diamond!e\leftrightarrow\phantom{!}e+\mbox{\ensuremath{!}}\Diamond e$}} & {\scriptsize{$!\square!e\leftrightarrow\mbox{!}e+\phantom{!}\square e$}} & {\scriptsize{$(\mathrm{V,D,W,C})$}} & {\scriptsize{$\varrho_{22}$}}\tabularnewline
{\scriptsize{$\varrho_{22}$}} & {\scriptsize{$!\Diamond!e\leftrightarrow\mbox{!}e+\phantom{!}\Diamond e$}} & {\scriptsize{$!\square!e\leftrightarrow\phantom{!}e+\mbox{!}\square e$}} & {\scriptsize{$(\mathrm{C,D,V,W})$}} & {\scriptsize{$\varrho_{21}$}}\tabularnewline
{\scriptsize{$\varrho_{23}$}} & {\scriptsize{$!\Diamond!e\leftrightarrow\mbox{!}e+\mbox{\ensuremath{!}}\Diamond e$}} & {\scriptsize{$!\square!e\leftrightarrow\phantom{!}e+\phantom{!}\square e$}} & {\scriptsize{$(\mathrm{V,W,C,D})$}} & {\scriptsize{$\varrho_{20}$}}\tabularnewline
\hline 
\noalign{\vskip1pt}
\end{tabular}
\par\end{center}

\vspace{1mm}

\begin{elabeling}{00.0000.0000.0000}
\item [{\textbf{\textsc{Theorem~11}}}] $\mathcal{R}_{p}(1)$\emph{ is
a group isomorphic to the symmetric group $\mathrm{S}_{4}$.}
\end{elabeling}
\textbf{\textsc{Proof.}}\qquad{}We build the composition table for  $\mathcal{R}_{p}(1)$,
which, although laborious, can be performed by hand.
By setting $p=e$, $q=\lozenge e$, $r=\lozenge!e$,
we basically verify Boolean equivalences in $p,q,r$,
as in the examples below:\vspace{-1mm}
\begin{eqnarray*}
\varrho_{2}\,\varrho_{3} & \approx & \eta_{2}(e,\lozenge e,\lozenge!e)*\varrho_{3}\enskip\approx\enskip\eta_{2}(e,\varrho_{3}(e),\varrho_{3}(!e))\\
 & \approx & \varrho_{3}(e)\leftrightarrow e+\varrho_{3}(!e)\\
 & \approx & (\lozenge e\leftrightarrow e+\mbox{!}\lozenge!e)\leftrightarrow e+(\lozenge!e\leftrightarrow\mbox{!}e+\mbox{!}\lozenge e)\\
 & \approx & (q\leftrightarrow p+\mbox{!}r)\leftrightarrow p+(r\leftrightarrow\mbox{!}p+\mbox{!}q)\enskip\approx\enskip q\leftrightarrow p\\
 & \approx & \lozenge e\leftrightarrow e\enskip\approx\enskip\varrho_{1}
\end{eqnarray*}
\begin{eqnarray*}
\varrho_{8}\,\varrho_{9} & \approx & \eta_{8}(e,\lozenge e,\lozenge!e)*\varrho_{9}\enskip\approx\enskip\eta_{8}(e,\varrho_{9}(e),\varrho_{9}(!e))\\
 & \approx & !\varrho_{9}(e)\leftrightarrow e+\varrho_{9}(!e)\\
 & \approx & \mbox{\ensuremath{!}}(!\lozenge e\leftrightarrow e+\mbox{!}\lozenge!e)\leftrightarrow e+(!\lozenge!e\leftrightarrow\mbox{\mbox{!}}e+\mbox{!}\lozenge e)\\
 & \approx & \mbox{\ensuremath{!}}(!q\leftrightarrow p+\mbox{!}r)\leftrightarrow p+(!r\leftrightarrow\mbox{\mbox{!}}p+\mbox{!}q)\enskip\approx\enskip q\\
 & \approx & \lozenge e\enskip\approx\enskip\varrho_{0}
\end{eqnarray*}\vspace{-8mm}

Next, we need the calculations from Table 3.\vspace{-2mm}

\begin{center}
\begin{tabular}{|c|c|c|c|c|}
\multicolumn{5}{c}{{\footnotesize{Table 3: Products of UR formulas with complemental
arguments}}}\tabularnewline[2mm]
\hline 
{\scriptsize{ $\varrho_{i}$}} & {\scriptsize{$\varrho_{i}(e)\varrho_{i}(!e)$}} & {\scriptsize{$\varrho_{i}(e)!\varrho_{i}(!e)$}} & {\scriptsize{$!\varrho_{i}(e)\varrho_{i}(!e)$}} & {\scriptsize{$!\varrho_{i}(e)!\varrho_{i}(!e)$}}\tabularnewline
\hline 
\noalign{\vskip1pt}
{\scriptsize{$\varrho_{0}$}} & {\scriptsize{$\Diamond e\Diamond!e$}} & {\scriptsize{$\Diamond e!\Diamond!e$}} & {\scriptsize{$!\Diamond e\Diamond!e$}} & {\scriptsize{$!\Diamond e!\Diamond!e$}}\tabularnewline
{\scriptsize{$\varrho_{1}$}} & {\scriptsize{$e\phantom{!}\Diamond e!\Diamond!e+\mbox{!}e!\Diamond e\phantom{!}\Diamond!e$}} & {\scriptsize{$e\phantom{!}\Diamond e\phantom{!}\Diamond!e+\mbox{!}e!\Diamond e!\Diamond!e$}} & {\scriptsize{$e!\Diamond e!\Diamond!e+\mbox{!}e\phantom{!}\Diamond e\phantom{!}\Diamond!e$}} & {\scriptsize{$e!\Diamond e\phantom{!}\Diamond!e+\mbox{!}e\phantom{!}\Diamond e!\Diamond!e$}}\tabularnewline
{\scriptsize{$\varrho_{2}$}} & {\scriptsize{$\Diamond e\Diamond!e$}} & {\scriptsize{$e\phantom{!}\Diamond e!\Diamond!e+\mbox{!}e!\Diamond e!\Diamond!e$}} & {\scriptsize{$e!\Diamond e!\Diamond!e+\mbox{!}e!\Diamond e\phantom{!}\Diamond!e$}} & {\scriptsize{$e!\Diamond e\phantom{!}\Diamond!e+\mbox{!}e\phantom{!}\Diamond e!\Diamond!e$}}\tabularnewline
{\scriptsize{$\varrho_{3}$}} & {\scriptsize{$e\phantom{!}\Diamond e!\Diamond!e+\mbox{!}e!\Diamond e\phantom{!}\Diamond!e$}} & {\scriptsize{$e\phantom{!}\Diamond e\phantom{!}\Diamond!e+\mbox{!}e\phantom{!}\Diamond e!\Diamond!e$}} & {\scriptsize{$e!\Diamond e\phantom{!}\Diamond!e+\mbox{!}e\phantom{!}\Diamond e\phantom{!}\Diamond!e$}} & {\scriptsize{$!\Diamond e!\Diamond!e$}}\tabularnewline
{\scriptsize{$\varrho_{4}$}} & {\scriptsize{$\Diamond e\Diamond!e$}} & {\scriptsize{$e!\Diamond e!\Diamond!e+\mbox{!}e\phantom{!}\Diamond e!\Diamond!e$}} & {\scriptsize{$e!\Diamond e\phantom{!}\Diamond!e+\mbox{!}e!\Diamond e!\Diamond!e$}} & {\scriptsize{$e\phantom{!}\Diamond e!\Diamond!e+\mbox{!}e!\Diamond e\phantom{!}\Diamond!e$}}\tabularnewline
{\scriptsize{$\varrho_{5}$}} & {\scriptsize{$e!\Diamond e\phantom{!}\Diamond!e+\mbox{!}e\phantom{!}\Diamond e!\Diamond!e$}} & {\scriptsize{$e\phantom{!}\Diamond e!\Diamond!e+\mbox{!}e\phantom{!}\Diamond e\phantom{!}\Diamond!e$}} & {\scriptsize{$e\phantom{!}\Diamond e\phantom{!}\Diamond!e+\mbox{!}e!\Diamond e\phantom{!}\Diamond!e$}} & {\scriptsize{$!\Diamond e!\Diamond!e$}}\tabularnewline
\hline 
\noalign{\vskip1pt}
{\scriptsize{$\varrho_{6}$}} & {\scriptsize{$!\Diamond e!\Diamond!e$}} & {\scriptsize{$!\Diamond e\Diamond!e$}} & {\scriptsize{$\Diamond e!\Diamond!e$}} & {\scriptsize{$\Diamond e\Diamond!e$}}\tabularnewline
{\scriptsize{$\varrho_{7}$}} & {\scriptsize{$e!\Diamond e\phantom{!}\Diamond!e+\mbox{!}e\phantom{!}\Diamond e!\Diamond!e$}} & {\scriptsize{$e!\Diamond e!\Diamond!e+\mbox{!}e\phantom{!}\Diamond e\phantom{!}\Diamond!e$}} & {\scriptsize{$e\phantom{!}\Diamond e\phantom{!}\Diamond!e+\mbox{!}e!\Diamond e!\Diamond!e$}} & {\scriptsize{$e\phantom{!}\Diamond e!\Diamond!e+\mbox{!}e!\Diamond e\phantom{!}\Diamond!e$}}\tabularnewline
{\scriptsize{$\varrho_{8}$}} & {\scriptsize{$e!\Diamond e\phantom{!}\Diamond!e+\mbox{!}e\phantom{!}\Diamond e!\Diamond!e$}} & {\scriptsize{$e!\Diamond e!\Diamond!e+\mbox{!}e!\Diamond e\phantom{!}\Diamond!e$}} & {\scriptsize{$e\phantom{!}\Diamond e!\Diamond!e+\mbox{!}e!\Diamond e!\Diamond!e$}} & {\scriptsize{$\Diamond e\Diamond!e$}}\tabularnewline
{\scriptsize{$\varrho_{9}$}} & {\scriptsize{$!\Diamond e!\Diamond!e$}} & {\scriptsize{$e!\Diamond e\phantom{!}\Diamond!e+\mbox{!}e\phantom{!}\Diamond e\phantom{!}\Diamond!e$}} & {\scriptsize{$e\phantom{!}\Diamond e\phantom{!}\Diamond!e+\mbox{!}e\phantom{!}\Diamond e!\Diamond!e$}} & {\scriptsize{$e\phantom{!}\Diamond e!\Diamond!e+\mbox{!}e!\Diamond e\phantom{!}\Diamond!e$}}\tabularnewline
{\scriptsize{$\varrho_{10}$}} & {\scriptsize{$e\phantom{!}\Diamond e!\Diamond!e+\mbox{!}e!\Diamond e\phantom{!}\Diamond!e$}} & {\scriptsize{$e!\Diamond e\phantom{!}\Diamond!e+\mbox{!}e!\Diamond e!\Diamond!e$}} & {\scriptsize{$e!\Diamond e!\Diamond!e+\mbox{!}e\phantom{!}\Diamond e!\Diamond!e$}} & {\scriptsize{$\Diamond e\Diamond!e$}}\tabularnewline
{\scriptsize{$\varrho_{11}$}} & {\scriptsize{$!\Diamond e!\Diamond!e$}} & {\scriptsize{$e\phantom{!}\Diamond e\phantom{!}\Diamond!e+\mbox{!}e!\Diamond e\phantom{!}\Diamond!e$}} & {\scriptsize{$e\phantom{!}\Diamond e!\Diamond!e+\mbox{!}e\phantom{!}\Diamond e\phantom{!}\Diamond!e$}} & {\scriptsize{$e!\Diamond e\phantom{!}\Diamond!e+\mbox{!}e\phantom{!}\Diamond e!\Diamond!e$}}\tabularnewline
\hline 
\noalign{\vskip1pt}
{\scriptsize{$\varrho_{12}$}} & {\scriptsize{$\Diamond e\Diamond!e$}} & {\scriptsize{$!\Diamond e\Diamond!e$}} & {\scriptsize{$\Diamond e!\Diamond!e$}} & {\scriptsize{$!\Diamond e!\Diamond!e$}}\tabularnewline
{\scriptsize{$\varrho_{13}$}} & {\scriptsize{$e!\Diamond e\phantom{!}\Diamond!e+\mbox{!}e\phantom{!}\Diamond e!\Diamond!e$}} & {\scriptsize{$e\phantom{!}\Diamond e\phantom{!}\Diamond!e+\mbox{!}e!\Diamond e!\Diamond!e$}} & {\scriptsize{$e!\Diamond e!\Diamond!e+\mbox{!}e\phantom{!}\Diamond e\phantom{!}\Diamond!e$}} & {\scriptsize{$e\phantom{!}\Diamond e!\Diamond!e+\mbox{!}e!\Diamond e\phantom{!}\Diamond!e$}}\tabularnewline
{\scriptsize{$\varrho_{14}$}} & {\scriptsize{$\Diamond e\Diamond!e$}} & {\scriptsize{$e!\Diamond e\phantom{!}\Diamond!e+\mbox{!}e!\Diamond e!\Diamond!e$}} & {\scriptsize{$e!\Diamond e!\Diamond!e+\mbox{!}e\phantom{!}\Diamond e!\Diamond!e$}} & {\scriptsize{$e\phantom{!}\Diamond e!\Diamond!e+\mbox{!}e!\Diamond e\phantom{!}\Diamond!e$}}\tabularnewline
{\scriptsize{$\varrho_{15}$}} & {\scriptsize{$e!\Diamond e\phantom{!}\Diamond!e+\mbox{!}e\phantom{!}\Diamond e!\Diamond!e$}} & {\scriptsize{$e\phantom{!}\Diamond e\phantom{!}\Diamond!e+\mbox{!}e!\Diamond e\phantom{!}\Diamond!e$}} & {\scriptsize{$e\phantom{!}\Diamond e!\Diamond!e+\mbox{!}e\phantom{!}\Diamond e\phantom{!}\Diamond!e$}} & {\scriptsize{$!\Diamond e!\Diamond!e$}}\tabularnewline
{\scriptsize{$\varrho_{16}$}} & {\scriptsize{$\Diamond e\Diamond!e$}} & {\scriptsize{$e!\Diamond e!\Diamond!e+\mbox{!}e!\Diamond e\phantom{!}\Diamond!e$}} & {\scriptsize{$e\phantom{!}\Diamond e!\Diamond!e+\mbox{!}e!\Diamond e!\Diamond!e$}} & {\scriptsize{$e!\Diamond e\phantom{!}\Diamond!e+\mbox{!}e\phantom{!}\Diamond e!\Diamond!e$}}\tabularnewline
{\scriptsize{$\varrho_{17}$}} & {\scriptsize{$e\phantom{!}\Diamond e!\Diamond!e+\mbox{!}e!\Diamond e\phantom{!}\Diamond!e$}} & {\scriptsize{$e!\Diamond e\phantom{!}\Diamond!e+\mbox{!}e\phantom{!}\Diamond e\phantom{!}\Diamond!e$}} & {\scriptsize{$e\phantom{!}\Diamond e\phantom{!}\Diamond!e+\mbox{!}e\phantom{!}\Diamond e!\Diamond!e$}} & {\scriptsize{$!\Diamond e!\Diamond!e$}}\tabularnewline
\hline 
\noalign{\vskip1pt}
{\scriptsize{$\varrho_{18}$}} & {\scriptsize{$!\Diamond e!\Diamond!e$}} & {\scriptsize{$\Diamond e!\Diamond!e$}} & {\scriptsize{$!\Diamond e\Diamond!e$}} & {\scriptsize{$\Diamond e\Diamond!e$}}\tabularnewline
{\scriptsize{$\varrho_{19}$}} & {\scriptsize{$e\phantom{!}\Diamond e!\Diamond!e+\mbox{!}e!\Diamond e\phantom{!}\Diamond!e$}} & {\scriptsize{$e!\Diamond e!\Diamond!e+\mbox{!}e\phantom{!}\Diamond e\phantom{!}\Diamond!e$}} & {\scriptsize{$e\phantom{!}\Diamond e\phantom{!}\Diamond!e+\mbox{!}e!\Diamond e!\Diamond!e$}} & {\scriptsize{$e!\Diamond e\phantom{!}\Diamond!e+\mbox{!}e\phantom{!}\Diamond e!\Diamond!e$}}\tabularnewline
{\scriptsize{$\varrho_{20}$}} & {\scriptsize{$e\phantom{!}\Diamond e!\Diamond!e+\mbox{!}e!\Diamond e\phantom{!}\Diamond!e$}} & {\scriptsize{$e!\Diamond e!\Diamond!e+\mbox{!}e\phantom{!}\Diamond e!\Diamond!e$}} & {\scriptsize{$e!\Diamond e\phantom{!}\Diamond!e+\mbox{!}e!\Diamond e!\Diamond!e$}} & {\scriptsize{$\Diamond e\Diamond!e$}}\tabularnewline
{\scriptsize{$\varrho_{21}$}} & {\scriptsize{$!\Diamond e!\Diamond!e$}} & {\scriptsize{$e\phantom{!}\Diamond e!\Diamond!e+\mbox{!}e\phantom{!}\Diamond e\phantom{!}\Diamond!e$}} & {\scriptsize{$e\phantom{!}\Diamond e\phantom{!}\Diamond!e+\mbox{!}e!\Diamond e\phantom{!}\Diamond!e$}} & {\scriptsize{$e!\Diamond e\phantom{!}\Diamond!e+\mbox{!}e\phantom{!}\Diamond e!\Diamond!e$}}\tabularnewline
{\scriptsize{$\varrho_{22}$}} & {\scriptsize{$e!\Diamond e\phantom{!}\Diamond!e+\mbox{!}e\phantom{!}\Diamond e!\Diamond!e$}} & {\scriptsize{$e\phantom{!}\Diamond e!\Diamond!e+\mbox{!}e!\Diamond e!\Diamond!e$}} & {\scriptsize{$e!\Diamond e!\Diamond!e+\mbox{!}e!\Diamond e\phantom{!}\Diamond!e$}} & {\scriptsize{$\Diamond e\Diamond!e$}}\tabularnewline
{\scriptsize{$\varrho_{23}$}} & {\scriptsize{$!\Diamond e!\Diamond!e$}} & {\scriptsize{$e\phantom{!}\Diamond e\phantom{!}\Diamond!e+\mbox{!}e\phantom{!}\Diamond e!\Diamond!e$}} & {\scriptsize{$e!\Diamond e\phantom{!}\Diamond!e+\mbox{!}e\phantom{!}\Diamond e\phantom{!}\Diamond!e$}} & {\scriptsize{$e\phantom{!}\Diamond e!\Diamond!e+\mbox{!}e!\Diamond e\phantom{!}\Diamond!e$}}\tabularnewline
\hline 
\noalign{\vskip1pt}
\end{tabular}

\par\end{center}

These equivalences too can be verified with simple Boolean calculations.
We also have $\Diamond e\Diamond!e\,\approx\, e\Diamond e\Diamond!e+\mbox{!}e\Diamond e\Diamond!e$
and 3 other similar equivalences.
Then the pattern that we observe is a formula that we will use shortly:\vspace{-4mm}

\begin{equation}
\text{¡}\varrho(e)\,\text{¡}\varrho(!e)\,\approx\, e\,\text{¡}\Diamond e\,\text{¡}\Diamond!e+\mbox{!}e\,\text{¡}\Diamond e\,\text{¡}\Diamond!e\label{eq:MeM!e}
\end{equation}

\pagebreak{}

Although not reflected in the above formula, there are obvious correlations
between the particular occurrences of `$\text{¡}$' and they also
depend on $\varrho$. However, we chose not to encumber the notation
with indexes on `$\text{¡}$'.
The reader only needs to keep in mind that here, as well as in subsequent formulas, the pseudo-operator `$\text{¡}$' instance typically depends on its position.

To prove that $\mathcal{R}_p(1)$ is isomorphic to $\mathrm{S}_4$
we consider its group action on the set
$\{\mathrm{W},\mathrm{D},\mathrm{C},\mathrm{V}\}$,
where $\mathrm{W}=\lozenge\mathit{1}\lozenge\mathit{0}$,
$\mathrm{D}=\lozenge\mathit{1}\mbox{!}\lozenge\mathit{0}$,
$\mathrm{C}=\mbox{!}\lozenge\mathit{1}\lozenge\mathit{0}$ and
$\mathrm{V}=\mbox{!}\lozenge\mathit{1}\mbox{!}\lozenge\mathit{0}$ are the axioms
of the 4 atoms of the lattice of $\mathbf{E}[0,1]$ systems from Figure 1.
We take $e=\mathit{1}$ in Table 3, which yields all the combinations
$(\text{¡}\lozenge\mathit{1}\text{¡}\lozenge\mathit{0})*\varrho_i$.
Then we see that $\mathcal{R}(1)$ generates precisely the 24 permutations of
$(\mathrm{W},\mathrm{D},\mathrm{C},\mathrm{V})$,
as per column 4 of Table 2.
\textsc{\scriptsize{\hfill{}$\blacksquare$}}\vspace{3mm}

We now show that the prime URs carry these 24 symmetries over from $\mathbf{E}[0,1]$
to all the other $\mathbf{E}[v,d]$ contexts.
We can actually define the state tuples $(\text{¡}\lozenge\mathit{1},\text{¡}\lozenge\mathit{0})$
as the \emph{primary labels} W, D, C and V of the minterms.
While prime substitutions permute minterms and preserve their labels,
we show that prime URs permute minterms and their labels as per Table 2.
\bgroup \renewcommand*\theenumi{\alph{enumi}}\renewcommand*\labelenumi{\theenumi)}

\begin{elabeling}{00.0000.0000.0000}
\item [{\textbf{\textsc{Theorem~12}}}] \noindent \emph{Let }$\varrho\in\mathcal{R}_{p}(1)$\emph{.
Then for every context }$\mathbf{E}[v,d]$:\end{elabeling}\vspace{-6mm}
\begin{enumerate}
\item $f_{\varrho}:[\mathit{1}]\rightarrow[\mathit{1}]$
\emph{ with }$f_{\varrho}(\mu)=\mu*\varrho$\emph{ for any minterm }$\mu$
\emph{is a bijection on the set of minterms.}\vspace{-2mm}
\item $f_{\varrho}:\Omega([\mathit{1}])\rightarrow\Omega([\mathit{1}]$)\emph{ with }$f_{\varrho}(\omega)=\omega*\varrho$\emph{ for any prime orbit }$\omega$
\emph{is a bijection on
the set of prime orbits.}
\end{enumerate}\vspace{-1mm}
\textbf{\textsc{Proof.}}\qquad{}For a), we first show by induction
on the modal level $d$ that $\varrho$ transforms any minterm $\mu$
into a single minterm $\mu_\varrho=\mu\,*\,\varrho$.

If $d=0$ the claim is
obvious, since level 0 minterms are unchanged by URs. Assuming the claim
holds up to some level $d$, it follows immediately that $f_{\varrho}$ is a
bijection on $\mathbf{E}[v,d]$. Let $\mu=\epsilon\prod_{l}\text{¡}\Diamond\phi_{l}$
be a level $d+1$ minterm, where $\epsilon$ is a level 0 minterm and
$\phi_{l}\in\mathbf{E}[v,d]$ for $1\leq l\leq|\mathbf{E}[v,d]|$.
Then:\egroup\vspace{-2mm}

\[
\mu*\varrho\approx\epsilon\prod_{l}\text{¡}(\Diamond\phi_{l})*\varrho\approx\epsilon\prod_{l}\text{¡}\varrho(\phi_{l}*\varrho)
\]\vspace{-3mm}

By grouping the modal factors
in pairs that have complemental formulas $\phi_{l}$ under the modal
operators, while at the same time adjusting the range of $l$ (to half
the initial range), we get:\vspace{-2mm}

\[
\mu*\varrho\approx\epsilon\prod_{l}(\,\text{¡}\varrho(\phi_{l}*\varrho)\,\text{¡}\varrho(!\phi_{l}*\varrho)\,)
\]\vspace{-3mm}

In the above we can use (\ref{eq:MeM!e}) with $\phi_{l}*\varrho$ as $e$
to write:\vspace{-5mm}

\[
\mu*\varrho\approx\epsilon\prod_{l}(\,(\phi_{l}*\varrho)\text{¡}\Diamond(\phi_{l}*\varrho)\,\text{¡}\Diamond(!\phi_{l}*\varrho)+\mbox{!}(\phi_{l}*\varrho)\text{¡}\Diamond(\phi_{l}*\varrho)\,\text{¡}\Diamond(!\phi_{l}*\varrho)\,)
\]

But the modal degree of $\phi_{l}$ is $\leq d$, therefore
\[
\{\Diamond(\phi_{l}*\varrho):\phi_{l}\in\mathbf{E}[v,d]\}=
\{\Diamond\phi_{l}:\phi_{l}\in\mathbf{E}[v,d]\}
\]
is the set of all the level $d+1$ modal factors,
since the level $d$ formulas $\phi_{l}$ are permuted by $\varrho$.
Then we can renumber all $\phi_{l}$ appropriately such that:\vspace{-1mm}

\begin{equation}
\mu*\varrho\approx\epsilon\prod_{l}(\phi_{l}\text{¡}\Diamond\phi_{l}\,\text{¡}\Diamond!\phi_{l}+\mbox{!}\phi_{l}\text{¡}\Diamond\phi_{l}\,\text{¡}\Diamond!\phi_{l})\approx\epsilon\prod_{l}\psi_{l}(\phi_{l})\label{eq:OneMinterm}
\end{equation}
\vspace{-1mm}

To prove that this formula represents a single minterm, observe that
the modal factors on level $d+1$ actually include all the modal factors
from levels 1 to $d$ (because some level $d+1$ formulas reduce
to level $d$ etc.). So in (\ref{eq:OneMinterm}) we start the multiplication
with $\epsilon$, which is a single level 0 minterm, and we consider
any level 1 factor $\psi_{l}$ from the product. Its sub-formulas
$\phi_{l}$ and $!\phi_{l}$ have modal degree 0, but only one of
these two minmatrices can include $\epsilon$, hence $\epsilon$ reduces
$\psi_{l}$ to (the corresponding) $\text{¡}\Diamond\phi_{l}\,\text{¡}\Diamond!\phi_{l}$.
Since (\ref{eq:OneMinterm}) contains every level 1 modal factor in
some $\psi_{l}$, it follows that the product of $\epsilon$ and all
level 1 factors $\psi_{l}$ is in fact a single level 1 minterm $\epsilon'$.
Then we can repeat this argument for $\epsilon'$ and all level 2
factors $\psi_{l}$ to infer that their product is a single level
2 minterm $\epsilon''$, and so on up to level $d+1$ where we conclude
that the whole formula (\ref{eq:OneMinterm}) is indeed a single level
$d+1$ minterm $\mu_\varrho$.

Finally, since $\varrho$ is invertible on the finite set of minterms
from $\mathbf{E}[v,d]$, it is a bijection.

For b), assume that minterms $\mu_{1}$ and $\mu_{2}$ are included
in a prime orbit $\omega$. Then there is a prime substitution $\varsigma_{i}$
such that $\mu_{1}\circ\varsigma_{i}\approx\mu_{2}$, therefore by
Lemma 7 $(\mu_{1}\circ\varsigma_{i})*\varrho\approx\mu_{2}*\varrho$
and by Lemma 8 $(\mu_{1}*\varrho)\circ(\varsigma_{i}*\varrho)\approx\mu_{2}*\varrho$.
But since $\varsigma_{i}$ is a level 0 substitution $\varsigma_{i}*\varrho\approx\varsigma_{i}$,
hence $(\mu_{1}*\varrho)\circ\varsigma_{i}\approx\mu_{2}*\varrho$,
which implies that $\mu_{1}*\varrho$ and $\mu_{2}*\varrho$ are both
in $\omega_{\varrho}\triangleq\omega*\varrho$.

A similar argument for $\varrho^{-1}$ and $\varsigma_{i}^{-1}$ combined
with claim a) above proves that prime orbits $\omega$ and $\omega_{\varrho}$
correspond through $\varrho$ and $\varrho^{-1}$ and have the same
number of minterms. Since $\varrho$ is now invertible on the finite
set of prime orbits from $\mathbf{E}[v,d]$, it is a bijection.\textsc{\scriptsize{\hfill{}$\blacksquare$}}\vspace{2mm}

\begin{elabeling}{00.0000.0000.0000}
\item [{\textbf{\textsc{Corollary~13}}}] \noindent $f_{\varrho}:\mathbf{E}[v,d]\rightarrow\mathbf{E}[v,d]$\emph{
defined as }$f_{\varrho}(\varphi)=\varphi*\varrho$\emph{ is a lattice
automorphism.}\vspace{1mm}

\end{elabeling}

As a side note, we mention that for a non-prime UR $\rho$, the function
$f_{\rho}:\Omega([\mathit{1}])\rightarrow\wp(\Omega([\mathit{1}]$))
defined as $f_{\rho}(\omega)=\omega*\rho$ transforms any prime orbit
$\omega$ into a (possibly empty) set of complete prime orbits $f_{\rho}(\omega)$,
such that if $\omega\not=\omega'$ then $f_{\rho}(\omega)\cap f_{\rho}(\omega')=\varnothing$
and $\bigcup_{\omega\in\Omega([\mathit{1}])}f_{\rho}(\omega)=\Omega([\mathit{1}])$.
But since we are not going to make use of this result, we leave the
proof to the reader.
\begin{elabeling}{00.0000.0000.0000}

\pagebreak{}

\item [{\textbf{\textsc{Theorem~14}}}] \emph{Let }$\varrho\in\mathcal{R}_{p}(1)$\emph{.
Then for any context }$\mathbf{E}[v,d]$\emph{, the function }$f_{\varrho}:\mathbf{E}\cml v,d\cmr \rightarrow\mathbf{E}\cml v,d\cmr$\emph{
defined as }$f_{\varrho}(\xi)=\xi*\varrho$\emph{ is a lattice automorphism.}
\end{elabeling}
\textbf{\textsc{Proof.}}\qquad{}For $d=0$ the theorem applies trivially
because these contexts reduce to PC, where
only $[\mathit{0}]$ and $[\mathit{1}]$ are CMMs and URs leave all formulas unchanged.
So we need to consider only contexts with $d>0$.

Basically, we must show that any prime UR $\varrho$ transforms a CMM from $\mathbf{E}\cml v,d\cmr$
into another CMM from $\mathbf{E}\cml v,d\cmr$, that the transformation
is invertible and that it preserves the lattice operations that we defined
for $\mathbf{E}\cml v,d\cmr$.

By Theorem 12, $\varrho$ is a bijection on the prime orbits,
so it transforms any CMM $\xi$ into a minmatrix $\xi*\varrho$ with the same number of
complete prime orbits. But since $\xi*\varrho$ may, in principle,
collapse, we need to prove that it is indeed another CMM. For this,
we show that in general a minmatrix $_{v}^{d}[\varphi]$ collapses iff $_{v}^{d}[\varphi*\varrho]$
collapses in a similar way, namely to prime orbits that correspond
through $\varrho$ and $\varrho^{-1}$.

Fix a context $\mathbf{E}[v,d]$ and consider an axiom $\mathrm{S}$
together with its associated $\mathbf{E}[v,d]$ system $\mathbf{S}$,
as well as their correspondents $\mathrm{S}_{\varrho}\triangleq\mathrm{S}*\varrho$
and $\mathbf{S_{\varrho}}$. Let $_{v}^{d}\varphi$ be a theorem of
$\mathbf{S}$, then consider any formal proof of $\varphi$, denoted
as the \emph{left hand proof}. Using the definition of $\varrho$,
syntactically transform every line of this proof into its corresponding
formula through $\varrho$ and denote the resulting sequence of formulas
as the \emph{right-hand proof}. We claim that the latter is a valid
proof of $_{v}^{d}\varphi_{\varrho}\triangleq\,_{v}^{d}\varphi*\varrho$
in $\mathbf{S}_{\mathbf{\varrho}}$ (although strictly speaking it
needs to be augmented with the formal proof of the arguments 1--3
below, every time they are applied). This is to say that the $\varrho$-transformed
formulas in the right-hand proof remain compatible with the same inference
rules of $\mathbf{E}$ that were applied to the corresponding lines
in the left-hand proof. But the inference rules of $\mathbf{E}$ are
as follows:\vspace{-1mm}
\begin{enumerate}
\item MP:
If the left hand proof contains $\alpha$, $\alpha\rightarrow\beta$
and $\beta$, then the right-hand proof contains $\alpha*\varrho$,
$(\alpha\rightarrow\beta)*\varrho$ and $\beta*\varrho$. But by UR-7
we have $(\alpha\rightarrow\beta)*\varrho\approx\alpha*\varrho\rightarrow\beta*\varrho$,
so the detachment of $\beta*\varrho$ is justified.
\item US:
If the left-hand proof contains $\alpha$ and $\alpha\circ\sigma$
then the right hand proof contains $\alpha*\varrho$ and $(\alpha\circ\sigma)*\varrho$.
But by Lemma 8 the latter is $(\alpha*\varrho)\circ(\sigma*\varrho)$, which
is a valid application of US in $\alpha*\varrho$.
\item RE-$\Diamond$:
If the left-hand proof contains $\alpha\leftrightarrow\beta$
and $\Diamond\alpha\leftrightarrow\Diamond\beta$ then the right hand
proof contains $(\alpha\leftrightarrow\beta)*\varrho$ and $(\Diamond\alpha\leftrightarrow\Diamond\beta)*\varrho$.
But the former is $\alpha*\varrho\leftrightarrow\beta*\varrho$ and
the latter is $\varrho(\alpha*\varrho)\leftrightarrow\varrho(\beta*\varrho)$,
which is a valid application of EQ (and implicitly RE)
on the right-hand side.
\end{enumerate}\vspace{-1mm}
If we assume that the right-hand proof is valid then the left-hand
proof is its $\varrho^{-1}$-transform, so by a similar argument it
is valid.

Back to $\xi=\cml\mathbf{S}\cmr$,
as a theorem, $\xi$ has a proof in $\mathbf{S}$.
The corresponding $\xi*\varrho$ has the same number of minterms and prime orbits as $\xi$,
and its proof in $\mathbf{S_\mathbf{\varrho}}$
shows that $\cml\mathbf{S_\mathbf{\varrho}}\cmr\subseteq\xi*\varrho$.
But if $\cml\mathbf{S_\mathbf{\varrho}}\cmr\subsetneq\xi*\varrho$,
then correspondingly
$\cml\mathbf{S_\mathbf{\varrho}}\cmr*\varrho^{-1}\subsetneq\xi$ and $\xi$ would collapse by intersection with theorem
$\cml\mathbf{S_\mathbf{\varrho}}\cmr*\varrho^{-1}$,
against our assumption that it is a CMM.
Thus, $\xi=\cml\mathbf{S}\cmr$ iff
$\xi*\varrho=\cml\mathbf{S_{\varrho}}\cmr$,
both consisting of prime orbits that correspond through $\varrho$ and
$\varrho^{-1}$. 

From Lemma 7, $\cml\mathbf{F}\cmr*\varrho=\cml\mathbf{F}\cmr$
and $\cml\mathbf{E}\cmr*\varrho=\cml\mathbf{E}\cmr$.
Next, consider the $\mathbf{E}[v,d]$ systems $\mathbf{S}$, $\mathbf{S'}$, $\mathbf{S''}$ with $\mathbf{\cml S\cmr}=\mathbf{\cml S'\cmr}\wedge\mathbf{\cml S''\cmr}$
on the left-hand side and their $\varrho$-correspondents on the right-hand side.
We have $\cml\mathbf{S}\cmr\subseteq\cml\mathbf{S'}\cmr\cap\cml\mathbf{S''}\cmr$,
where the intersection may collapse to some CMM as per Theorem 2 b).
However, every minmatrix that participates in the intersection that defines
the left-hand CMM $\cml\mathbf{S}\cmr$
(and whose proof may now include theorems from both $\mathbf{S'}$ and $\mathbf{S''}$)
has a right-hand correspondent,
therefore $\cml\mathbf{S}\cmr$ and $\cml\mathbf{S_\mathbf{\varrho}}\cmr$,
which must exist due to the finiteness of the context,
also correspond.

The dual argument applies to the $\vee$ operator, where we actually
have $\cml\mathbf{S}\cmr=\cml\mathbf{S'}\cmr\cup\cml\mathbf{S''}\cmr$,
because Theorem 2 a) shows that the union cannot collapse. Hence $\varrho$
preserves the lattice operations that we defined for $\mathbf{E}\cml v,d\cmr$,
which completes the proof\textsc{\scriptsize{.\hfill{}$\blacksquare$}}\vspace{3mm}

Observe that the gist of the above theorem is the left-hand-right-hand
correspondence between a proof and its $\varrho$-transform, which
we abbreviate as the \emph{LH-RH argument}. The intermediate formulas
in the proofs on both sides do not necessarily belong to the same
context as $\varphi$ and $\varphi_{\varrho}$, but since the UR properties
hold in every context, this has no impact on the conclusion.

Also observe that these automorphisms apply to base $\mathbf{E}$,
which has no additional axioms. But the LH-RH argument does not apply
to another base like, say, $\mathbf{K}$; for in that case a system
$\mathbf{S}_{\varrho}$ lacks the axiom $\mathrm{K}_{\varrho}$.

\section{Automorphisms of CExtE}

We now turn our attention to $\mathrm{CExt}\mathbf{E}$, which includes
all classical modal logics, whether or not finitely-axiomatizable.

As mentioned at the end of Section 2, all systems have a CMM in every
context. For any system $\mathbf{S}$, let $\mathcal{C}(\mathbf{S})=\{_{v}^{d}\cml\mathbf{S}\cmr:v,d\text{ integers, }v\geq0,d\geq0\}$.
(For the purpose of this set the CMMs from different contexts are
considered distinct elements even when equiprovable, since they have
a different DCF representation in every context. With this convention
we can avoid using tuples $(v,d,{}_{v}^{d}\cml S\cmr)$ instead.)
Then $\mathcal{C}(\mathbf{S})$ uniquely determines $\mathbf{S}$,
since any $\mathbf{S}$-theorem from a given context can be derived
from the context's CMM by PC-monotony. Hence $\mathbf{S}=\mathbf{S'}$
as sets of formulas iff $\mathcal{C}(\mathbf{S})=\mathcal{C}(\mathbf{S'})$.

We note that there are some constraints on the CMMs from $\mathcal{C}(\mathbf{S})$.
One cannot just randomly pick one CMM from every context, because
they all act as theorems and their combination might cause collapses.
For example, since the promoted $_{v}^{d}\cml\mathbf{S}\cmr$ is also a theorem
in $\mathbf{E}[v,d+1]$,
it always needs to include $_{v}^{d+1}\cml\mathbf{S}\cmr$,
otherwise the latter would collapse by their intersection.
Similarly, $_{v}^{d+1}\cml\mathbf{S}\cmr$ cannot be arbitrarily strong,
otherwise it may be used to collapse $_{v}^{d}\cml\mathbf{S}\cmr$.
But for the purpose of the next theorem, the exact nature of these
CMM constraints is not relevant. All that matters is that none of
the CMMs in $\mathcal{C}(\mathbf{S})$ collapses, as per their definition.
We call this the \emph{CMM compatibility} in the set $\mathcal{C}(\mathbf{S})$. 

For an arbitrary system $\mathbf{S}$ and $\varrho\in\mathcal{R}_{p}(1)$,
define $\mathbf{S}_{\varrho}=\{\varphi*\varrho:\varphi\in\mathbf{S}\}$
and $\mathcal{C}_{\varrho}(\mathbf{S})=\{_{v}^{d}\cml \mathbf{S}\cmr*\varrho:v,d\text{ integers, }v\geq0,d\geq0\}$.\vspace{-2mm} \bgroup \renewcommand*\theenumi{\alph{enumi}}\renewcommand*\labelenumi{\theenumi)}\vspace{2mm}
\begin{elabeling}{00.0000.0000.00000}
\item [{\textbf{\textsc{Theorem~15}}}] \noindent \emph{Let }$\varrho\in\mathcal{R}_{p}(1)$\emph{.
Then:}\end{elabeling}\vspace{-5mm}
\begin{enumerate}
\item \noindent $f_{\varrho}:\mathrm{CExt}\mathbf{E}\rightarrow\mathrm{CExt}\mathbf{E}$\emph{
defined as }$f_{\varrho}(\mathbf{S})=\mathbf{S}_{\varrho}$\emph{
is a lattice automorphism.}\vspace{-1mm}
\item \emph{The set }$\mathcal{A}_{p}(1)=\{f_{\varrho}:\mathcal{\varrho\in R}_{p}(1)\}$\emph{
is a group of automorphisms of }$\mathrm{CExt}\mathbf{E}$\emph{ that
is isomorphic to the symmetric group $\mathrm{S_{4}}$.}
\end{enumerate}\vspace{-2mm}
\textbf{\textsc{Proof.}}\qquad{}For a), we first show that $f_{\varrho}$
is well-defined. From Lemma 10 we see that all $\varphi$ and $\varphi_{\varrho}$
from $\mathbf{S}$ and $\mathbf{S}_{\varrho}$ respectively correspond
through $\varrho$ and $\varrho^{-1}$, hence so do the elements from
$\mathcal{C}(\mathbf{S})$ and $\mathcal{C}_{\varrho}(\mathbf{S})$.
We use the LH-RH argument to see that $\mathcal{C}_{\varrho}(\mathbf{S})$
is a set of compatible CMMs: if there were a proof that could collapse
a minmatrix $\xi *\varrho\in\mathcal{C}_{\varrho}(\mathbf{S})$, then by applying $\varrho^{-1}$
to the proof we would collapse the CMM $\xi\in\mathcal{C}(\mathbf{S})$,
contradicting its definition. Next, for any context, the remaining
theorems of $\mathbf{S}$ and $\mathbf{S_{\varrho}}$ are precisely
those minmatrices that include the CMMs from their context,
and they also correspond minterm by minterm as per Theorem 12.
Consequently $\mathbf{S_{\varrho}}$ is a valid system of $\mathrm{CExt}\mathbf{E}$,
uniquely determined by $\mathcal{C}(\mathbf{S_{\varrho}})=\mathcal{C}_{\varrho}(\mathbf{S})$.

By a similar reasoning on the correspondence between the sets $\mathcal{C}(\mathbf{S})$
and $\mathcal{C}(\mathbf{S_{\mathbf{\varrho}}})$ the function $f_{\varrho}$
is both injective and surjective.\egroup

Obviously $f_\varrho(\mathbf{F})=\mathbf{F}$
and $f_\varrho(\mathbf{E})=\mathbf{E}$.
To prove that $f_{\varrho}$ is a lattice automorphism, we note that
the lattice operations in $\mathrm{CExt}\mathbf{E}$, denoted as $\oplus$
and $\odot$, are \emph{not} the $\vee$ and $\wedge$ operations
that we defined for $\mathbf{E}\cml v,d\cmr$. Thus, $\mathbf{S}=\mathbf{S'}\oplus\mathbf{S''}$
is defined as taking the union $\mathbf{S'}\cup\mathbf{S''}$ and
performing the closure with respect to MP, US and RE-$\lozenge$.
But these are the very rules that make the LH-RH argument work.
In any given context both $\mathbf{\cml S'}\cmr$
and $\mathbf{\cml S''}\cmr$ are theorems of $\mathbf{S'}\cup\mathbf{S''}$,
hence so is $\mathbf{\cml S'}\cmr\cap\mathbf{\cml S''}\cmr$.
This either is the CMM $\mathbf{\cml S}\cmr=\mathbf{\cml S'}\oplus\mathbf{S''}\cmr$ or, by producing more theorems using the closure operation,
it collapses to the CMM.
In any case, by the LH-RH argument the $\varrho$-transformed theorems produce $\mathbf{\cml S_\mathbf{\varrho}}\cmr=\mathbf{\cml S_\mathbf{\varrho}'}\oplus\mathbf{S_\mathbf{\varrho}''}\cmr$.
Since this holds in any context, $\mathcal{C}(\mathbf{S})$
and $\mathcal{C}(\mathbf{S_{\mathbf{\varrho}}})$ must correspond,
hence $f_\varrho(\mathbf{S})=\mathbf{S_\mathbf{\varrho}}$.

Similar considerations apply to $\mathbf{S}=\mathbf{S'}\odot\mathbf{S''}$,
defined as $\mathbf{S'}\cap\mathbf{S''}$. In this case every theorem
must include both CMMs from the context, hence also $\mathbf{\cml S'}\cmr\cup\cml\mathbf{S''}\cmr$,
which cannot collapse, therefore $\mathbf{\cml S}\cmr=\mathbf{\cml S'}\odot\mathbf{S''}\cmr\approx\mathbf{\cml S'}\cmr\cup\cml\mathbf{S''}\cmr$
and the $\varrho$-correspondence between $\mathcal{C}(\mathbf{S})$
and $\mathcal{C}(\mathbf{S_{\mathbf{\varrho}}})$ still holds.\vspace{1mm}

For b), we observe that all $f_\varrho$ with $\varrho\in\mathcal{R}_{p}(1)$ are distinct, since they cause distinct permutations of the $\mathbf{E}[0,1]$ systems from Figure 1.
By Lemma 10 we have $(\varphi*\rho)*\rho'\approx\varphi*(\rho\rho')$ for every $\varphi$,
hence for every CMM,
therefore $f_{\varrho}f_{\varrho'}=f_{\varrho\varrho'}$. Thus, composition
in $\mathcal{A}_{p}(1)$ is well-defined and it has the same table
as $\mathcal{R}_{p}(1)$, which is\emph{ }isomorphic to the group $\mathrm{S_{4}}$
by Theorem 11.
\textsc{\scriptsize{\hfill{}$\blacksquare$}}{\scriptsize \par}

\section{Final Remarks}

The usual method for studying a lattice $\mathrm{NExt}\mathbf{S}$
of extensions of a system $\mathbf{S}$
is to consider the lattice as a whole.
In this paper we have presented a complemental view
that can be obtained from the countable-contextualization of formulas.
This approach provides some insight into $\mathrm{CExt}\mathbf{E}$,
by showing how the lattice of finitely-axiomatizable systems
is a sort of ``fractal refinement'' of the lattices of context CMMs,
whereafter $\mathrm{CExt}\mathbf{E}$ is obtained
by performing a ``closure'' operation,
which adds the systems that are not finitely-axiomatizable.
We shall also apply this method to $\mathrm{NExt}\mathbf{K}$ in \cite{Soncodi}.

Several topics are for further study.
As working with canonical can be tedious,
exemplifying some properties requires computer-aided calculations,
hence we can only hint at them here. Some other properties can be proven symbolically,
but this will be the subject of a future paper.

A first question is whether or not the above are all the automorphisms
of $\mathrm{CExt}\mathbf{E}$. A brute-force search reveals that among
our 256 URs, only the 24 identified prime URs are invertible.
The search algorithm uses straightforward calculations with Boolean
formulas in 3 variables, in the style of the proof of Theorem 11.
Equivalence can be shown by converting formulas to their Boolean DNF,
which is a well-known procedure.

Next, one can similarly investigate the \emph{extended URs} (XURs) of the form
$\chi(e)=\eta(e,\lozenge\mathit{1},\lozenge e,\lozenge!e,\lozenge\mathit{0})$, in which
case there are significantly more, yet similar calculations with Boolean
formulas $\eta$ in 5 variables. The algorithm needs to verify
all the compositions of the  $2^{32}$ XURs, but it can be sped up by
checking only those XURs that transform each minterm into a single
minterm. Then it turns out that there is a larger group of invertible XURs,
which includes $\mathcal{R}_p(1)$ and has $24\cdot2^{4}\cdot3^{4}=31,104$ elements,
for example: 
\[
\chi(e)=(\lozenge\mathit{1}\rightarrow e+(\lozenge e\leftrightarrow(\lozenge\mathit{0}\rightarrow\lozenge!e)))\rightarrow\lozenge e(\lozenge\mathit{1}\rightarrow e)
\]
\[
\chi^{-1}(e)=(\lozenge\mathit{1}\rightarrow e+(\lozenge e\leftrightarrow\lozenge\mathit{0}+\lozenge!e))\rightarrow\lozenge e(\lozenge\mathit{1}\rightarrow e)
\]

It is yet to be determined what this group is and what ``magic'' works behind its formulas.
All the prime XURs $\chi(e)$ have precisely half of the 32 minterms from $\mathbf{E}[1,1]$,
but at this point we do not have a construction rule for them.
The theory from Section 4 would need to be generalized for XURs,
and there are only a few places where it depends on the form of $\eta$.
Obviously these prime XURs cannot generate distinct automorphisms in all contexts, as for example the lattice from Figure 1 admits only 24 distinct ones. Still they may reveal additional complex symmetries of $\mathrm{CExt}\mathbf{E}$.

However, what we have defined so far can be considered
level 1 URs, namely those that leave level 0 formulas unchanged. The
question is open whether higher level URs can be defined. For example,
level 2 URs would be translations that leave level 0 and level 1 formulas
unchanged etc. It may also be possible to define non-context-preserving
URs, in which case one could find interesting $\mathrm{CExt}\mathbf{E}$
homomorphisms.

As a consequence of these automorphisms, it can be shown that $\mathrm{CExt}\mathbf{E}$
has precisely 4 co-atoms, namely the systems determined by the axioms
$\mathrm{W_w}=\lozenge p$, $\mathrm{D_d}=\lozenge p \leftrightarrow p$,
$\mathrm{C_c}=\lozenge p \leftrightarrow \mbox{!}p$ and $\mathrm{V_v}=\mbox{!}\lozenge p$.
This represents a generalization in $\mathbf{E}$ of Makinson's theorem mentioned in \cite{Chagrov2},
which states that $\mathrm{NExt}\mathbf{K}$ has only
2 co-atoms, $\mathbf{Triv}$ and $\mathbf{Ver}$.

For the normal modal system $\mathbf{K}$, using $\mathbf{E[}2,1]$
calculations it can be shown that all its non-trivial $\mathcal{R}_p(1)$
transforms $\mathbf{K}_{\varrho}$ are distinct and not normal. Then
this is a hint (but no proof) in support of the conjecture from \cite{Kracht}
that the lattice $\mathrm{NExt}\mathbf{K}$ is rigid.
Basically, in the $\mathrm{S}_{4}$-symmetrical $\mathrm{CExt}\mathbf{E}$, $\mathrm{NExt}\mathbf{K}$ is a small, ``oblique'' sublattice.

From Figure 1, which is essentially the diagram of a 4-dimensional hypercube,
we see that the prime URs correspond to the symmetries of the hypercube with a fixed vertex. These are known to be characterized by $\mathrm{S}_{4}$.
Regarding the transformed systems $\mathbf{K_{\mathbf{\varrho}}}$ mentioned above,
one can see that $\mathbf{K_{\mathbf{\varrho}_{18}}}$
can be described in terms of Kripke frames
where the semantics of the modal operators $\lozenge$ and $\square$
are swapped. Similarly, the semantics for $\mathbf{K_{\mathbf{\varrho}_{6}}}$
can be modified such that the valuation of $\Diamond p$ at some world
$w$ is $\mathit{1}$ iff $w$ sees at least one world where the valuation
of $p$ is $\mathit{0}$ (instead of $\mathit{1}$).
With such modified semantics, the transformed systems and their
extensions are characterized by the same Kripke frames
as their non-transformednormal counterparts.
It would be interesting to know what such modified Kripke
and also neighborhood semantics correspond to all $\varrho_{i}$
(or $\chi_i$). 

Lastly, we note that in general, when internal symmetries of an object
are revealed, numerous other symmetry-caused properties typically follow.
So in this respect we hope that our paper will inspire further investigations of
$\mathrm{CExt}\mathbf{E}$, potentially leading to new findings about its sublattice
$\mathrm{NExt}\mathbf{K}$.

\pagebreak{}

\raggedbottom


\AuthorAdressEmail{Adrian Soncodi}
{Lecturer, University of Texas, Dallas}
{acs151130@utdallas.edu\linebreak
soncodi@verizon.net}

\vspace{40pt}
\setlength{\parindent}{0ex}
\textbf{Springer acknowledgement notice:}
\\

This is a post-peer-review, pre-copyedit version of an article published in

\ \ \ \ \ \textit{Studia Logica} vol 104, November 2015.
\\

The final authenticated version is available online at:

\ \ \ \ \ https://dx.doi.org/10.1007/s11225-015-9638-8

\end{document}